\documentclass[11pt,twoside,reqno]{amsart}
\usepackage[all]{xy}   
\usepackage {amssymb,latexsym,amsthm,amsmath,setspace,bm}
\usepackage[colorlinks, linkcolor=blue, citecolor=red]{hyperref}
\usepackage{graphicx}
\usepackage{mathtools}
\usepackage{xcolor}
\usepackage[utf8]{inputenc}
\usepackage{tikz}
\usepackage{enumitem}

\newcommand{\secref}[1]{Section~\ref{#1}}
\newcommand{\thmref}[1]{Theorem~\ref{#1}}
\newcommand{\lemref}[1]{Lemma~\ref{#1}}
\newcommand{\remref}[1]{Remark~\ref{#1}}

\newcommand{\propref}[1]{Proposition~\ref{#1}}
\newcommand{\corref}[1]{Corollary~\ref{#1}}

\makeatletter
\def\imod#1{\allowbreak\mkern10mu({\operator@font mod}\,\,#1)}
\makeatother

\newtheorem{theorem}{Theorem}[subsection]
\newtheorem{lemma}[theorem]{Lemma}
\newtheorem{corollary}[theorem]{Corollary}
\newtheorem{proposition}[theorem]{Proposition}
\newtheorem*{theorem*}{Theorem}

\theoremstyle{definition}
\newtheorem{remark}[theorem]{Remark}

\newtheorem{definition}[theorem]{Definition}

\numberwithin{equation}{section}

\newcommand{\ignore}[1]{}

\newcommand{\nc}{\newcommand}

\nc{\Cal}{\cal} \nc{\Xp}[1]{X^+(#1)} \nc{\Xm}[1]{X^-(#1)}
\nc{\on}{\operatorname} \nc{\Z}{{\mathbb{Z}}} \nc{\J}{{\cal J}} \nc{\C}{{\mathbb{C}}} \nc{\Q}{{\bold Q}} \nc{\R}{{\mathbb{R}}}
\nc{\K}{\bold{\kappa}}

\nc{\F}{\operatorname{\boF}}
\nc{\pr}{\operatorname{\bold{pr}}}
\nc{\ev}{\operatorname{{ev}}}
\nc{\hook}{\lceil} 
\nc{\HLn}{\noindent\makebox[\linewidth]{\rule{\textwidth}{1pt}}}

\nc{\N}{{\Bbb N}} \nc\boa{\bold a} \nc\bob{\bold b} \nc\boc{\bold c} \nc\bod{\bold d} \nc\boe{\bold e} \nc\bof{\bold f} \nc\bog{\bold g}
\nc\boh{\bold h} \nc\boi{\bold i} \nc\boj{\bold j} \nc\bok{\bold k} \nc
\bol{\bold l} \nc\bom{\bold m} \nc\bon{\bold n} \nc\boo{\bold o}
\nc\bop{\bold p} \nc\boq{\bold q} \nc\bor{\bold r} \nc\bos{\bold s} \nc\boT{\bold t} \nc\boF{\bold F} \nc\bou{\bold u} \nc\bov{\bold v}
\nc\bow{\bold w} \nc\boz{\bold z} \nc\boy{\bold y} \nc\ba{\bold A} \nc\bb{\bold B} \nc\bc{\bold C} \nc\bd{\bold D} \nc\be{\bold E} \nc\bg{\bold
	G} \nc\bh{\bold H} \nc\bi{\bold I} \nc\bj{\bold J} \nc\bk{\bold K} \nc\bl{\bold L} \nc\bM{\bold M} \nc\bn{\bold N} \nc\bo{\bold O} \nc\bp{\bold
	P} \nc\bq{\bold Q} \nc\br{\bold R} \nc\bs{\bold S} \nc\bt{\bold T} \nc\bu{\bold U} \nc\bv{\bold V} \nc\bw{\bold W} \nc\bz{\bold Z} \nc\bx{\bold
	x} \nc\KR{\bold{KR}} \nc\rk{\bold{rk}} \nc\het{\text{ht }}

\nc\udi{\underline i} \nc\udj{\underline j}

\nc\fn{{fin}}  \nc\af{{aff}}  \nc\tr{{tor}} \nc\btilde{\bold{\tilde{\bold{H}}}}

\nc{\mpp}{\rotatebox[origin=c]{180}{\pm}}

\nc\eps{\epsilon}
\nc\toa{\tilde a} \nc\tob{\tilde b} \nc\toc{\tilde c} \nc\tod{\tilde d} \nc\toe{\tilde e} \nc\tof{\tilde f} \nc\tog{\tilde g} \nc\toh{\tilde h}
\nc\toi{\tilde i} \nc\toj{\tilde j} \nc\tok{\tilde k} \nc\tol{\tilde l} \nc\tom{\tilde m}  \nc\ton{\tilde n} \nc\too{\tilde o} \nc\toq{\tilde q}
\nc\tor{\tilde r} \nc\tos{\tilde s} \nc\toT{\tilde t} \nc\tou{\tilde u} \nc\tov{\tilde v} \nc\tow{\tilde w} \nc\toz{\tilde z}

\begin{document}
	\setcounter{section}{0}
	\setcounter{tocdepth}{1}
	
	\title{Demazure Filtrations of Tensor Product Modules and Character Formula}
	\author[Divya Setia]{Divya Setia}
	\author{Tanusree Khandai}
	\address{Indian Institute of Science Education and Research Mohali, Knowledge City,  Sector 81, S.A.S. Nagar 140306, Punjab, India}
	\email{divyasetia01@gmail.com}
	\email{tanusree@iisermohali.ac.in}
	\date{}
	\thanks{}


	
	\begin{abstract} 
		We study the structure of the finite-dimensional representations of $\mathfrak{sl}_2[t]$, the current Lie algebra type of $A_1$, which are obtained by taking tensor products of special Demazure modules. We show that these representations admit a Demazure flag and obtain a closed formula for the graded multiplicities of the level 2 Demazure modules in the filtration of the tensor product of two local Weyl modules for $\mathfrak{sl}_2[t]$. Furthermore, we derive an explicit expression for  graded character of the tensor product of a local Weyl module with an irreducible $\mathfrak{sl}_2[t]$ module. In conjunction with the results of \cite{MR3210603}, our findings provide evidence for the conjecture in \cite{9} that the tensor product of Demazure modules of levels m and n respectively has a filtration by Demazure modules of level m + n.
	\end{abstract}

	\maketitle
	\section{INTRODUCTION}
	
	Let $\mathfrak {g}$ be a finite-dimensional simple Lie algebra over the complex field $\mathbb C$ and $\mathfrak {g}[t]$ be its associated current algebra,  which is the Lie algebra of polynomimal mappings from $\mathbb C$ to $\mathfrak g$.
	Due to its connections with combinatorics, number theory and mathematical physics, the theory of finite-dimensional representations of current algebras has garnered significant attention in the  past two decades. Motivated by a conjecture in \cite{9} that suggests tensor products of Demazure modules with levels $m$ and $n$, respectively, have a filtration by Demazure modules of level $m+n$, in
	this paper we consider the class of finite-dimensional representations, namely those of $\mathfrak{sl}_2[t]$, which are derived from tensor products of certain Demazure modules of $\mathfrak{sl}_2[t]$ and prove results that provide evidence in support of the conjecture. 

	

	Let  $\hat{\mathfrak g}$ be an affine Kac-Moody Lie algebra and $\hat{\mathfrak b}$ be a standard positive Borel subalgebra of $\hat{\mathfrak g}$. Given an integrable highest weight irreducible representation $\mathcal V$ of $\hat{\mathfrak g}$, a Demazure module associated with $\mathcal V$, is defined as the $\hat{\mathfrak b}$-module generated by an extremal weight vector $u_0$ of $\mathcal V$. These modules are $\mathfrak g$-stable when the restriction of the weight of $u_0$ to a Cartan subalgebra of $\mathfrak g$ is anti-dominant and in this case the corresponding Demazure module is a module of the current algebra $\mathfrak g[t]$. A Demazure module is said to be of level $\ell$ if the central element of $\hat{\mathfrak g}$ acts on $\mathcal V$ by the scalar $\ell$. A $\hat{\mathfrak b}$-module is said to admit a Demazure flag if it has a filtration by submodules whose successive quotients are Demazure modules. While it is known that the tensor product of two Demazure modules does not in general have a Demazure flag, in \cite{MR2214249, MR2855081, MR1987017} it was proved that for simply-laced Kac–Moody Lie algebras,  the tensor product of a one-dimensional Demazure module with an arbitrary one admits a Demazure flag. Our study here extends this question to tensor products of level 1 Demazure modules with Demazure modules of arbitrary levels for current algebras of type $A_1$.

	The notion of Weyl modules for affine Kac-Moody Lie algebras was introduced in \cite{MR1850556}.  In the case when $\mathfrak g=\mathfrak{sl}_2$, it was proven that these modules are the classical  limit ($q\rightarrow 1$) of standard modules of the quantum affine algebras. Furthermore, it was shown that for a dominant integral weight $\lambda$, $W_{loc}(\lambda)$, the local Weyl module with highest weight $\lambda$,   is a finite dimensional graded $\mathfrak g[t]$-module generated by a highest weight vector of weight $\lambda$ and every finite-dimensional graded $\mathfrak g[t]$-module with highest weight $\lambda$ is a quotient of $W_{loc}(\lambda)$. Subsequently it was demonstrated (\cite{MR2323538, MR2271991, MR1953294, MR1771615}) that  for current algebra of type ADE, the local Weyl modules are in fact Demazure modules of level 1 and  their Demazure characters 
	coincide with  non-symmetric Macdonald polynomials, specialized at t = 0. 
	On the other hand, it was demonstrated in \cite{MR2855081} that the local Weyl modules for current algebras of non-simply laced Lie algebras $\mathfrak{g}$ have filtrations by Demazure modules whose multiplicity coincides with the multiplicity of certain Demazure modules in local Weyl modules of type $A_r$ for suitable $r$. 
	
	In \cite{MR3296163}, a  family of finite-dimensional quotients of the local Weyl modules, often referred to as Chari-Ventakesh Modules (in short CV-modules) was introduced. It was shown that these modules subsume many disparate classes of finite-dimensional graded representations of current algebras. By definition, a CV module $V(\xi)$ is associated to a family of partitions ${\bold{\xi}} = \{\xi_\alpha\}_{\alpha}$, which is indexed by the set of positive roots of $\mathfrak g$. It was proven in \cite[Theorem 2]{MR3296163} that the Demazure modules in various levels can be realized as CV-modules associated to a set of rectangular or special near rectangular partitions. The interpretation of Demazure modules as CV-module greatly simplified their defining relations and aided in the study of the structure of fusion product modules for $\mathfrak{sl}_2[t]$. In \cite{MR3210603} a necessary and sufficient condition was obtained for the existence of level $\ell$ Demazure flags in an arbitrary CV-module for $\mathfrak{sl}_{2}[t]$ and combinatorial descriptions of the graded multiplicities of Demazure modules in $V(\xi)$ was given. Ensuing research (\cite{MR3461134, MR4229660}) on the generating functions associated with the graded multiplicities of Demazure modules in special CV-modules established connections between these modules, number theory and combinatorics.  
	
	Motivated by the study initiated in \cite{9} on the structure of tensor product of two local Weyl modules for $\mathfrak{sl}_2[t]$, we begin by obtaining a presentation of the representations of $\mathfrak{sl}_{2}[t]$ that arise from 
	taking the tensor product of a local Weyl module with a CV module $V(\xi)$.  Specializing the partition $\xi$, we then explore the tensor product of local Weyl modules with irreducible $\mathfrak{sl}_2[t]$-modules and local Weyl modules for $\mathfrak{sl}_2[t]$. In each of these cases, we observe that the corresponding tensor product possesses a filtration by CV-modules. Using \cite[Theorem 3.3]{MR3210603}, we thus conclude that each such module has a filtration by Demazure modules at appropriate levels.
	
	For a given positive integer $m$, let $W_{loc}(m)$ denote the local Weyl module of $\mathfrak{sl}_2[t]$ with the highest weight $m$. Through explicit construction, we demonstrate that when $\xi=(n)$, the module $W_{loc}(m)\otimes V(\xi)$ admits a filtration by CV-modules with hook-type partitions. Utilizing this filtration, we derive expressions for the outer multiplicities in the tensor product module $W_{loc}(m\omega)\otimes ev_{0}V(n\omega)$. As a consequence we are able to 
	express the product of a Schur polynomial with a specialized Macdonald polynomial in terms of Schur polynomials.

	Finally we consider the tensor product of a local Weyl module with certain special Demazure modules. In \cite{MR3407180}, it had been shown that the truncated local Weyl module of $\mathfrak{sl}_2[t]$ of highest weight $r\omega$ can be realized as a CV-module.  Interestingly, we observe that the tensor product of two local Weyl modules have a filtration by truncated local Weyl modules of suitable weights. We determine the graded character of such truncated local Weyl modules in terms of level 2 Demazure modules and give the graded multicilities of level 2 Demazure modules in $W_{loc}(m)\otimes W_{loc}(n)$.
	
	The character formulas for certain tensor product modules  had been obtained in \cite{9}. Using an uniform approach we give a direct proof for the character formula of the modules which  we consider.  Then, using the Chari-Pressley-Loktev basis (\cite[Section 6]{MR1850556}, \cite[Theorem 2.1.3]{MR2271991}) and Chari-Venkatesh basis (\cite[Theorem 5]{MR3296163}) for local Weyl modules for $\mathfrak{sl}_2[t]$, the presentation of the tensor products modules we obtained and the character formulas obtained, we show that in the cases under consideration, the tensor product modules have a Demazure filtration of appropriate level. 

	The paper is structured as follows: In Section 2, we set the notations
	and recall the basic definitions and results essential for our paper. Section 3 begins with a recap of the definition of CV modules. We then proceed to give a presentation via generators and relations for the tensor product modules
	$W_{loc}(m) \otimes V(\xi)$. Moving on to Section 4, we systematically construct a decreasing chain of submodules of $W_{loc}(m\omega)\otimes ev_{0}V(n\omega)$ such that the successive quotients in this chain are isomorphic to $V(\xi)$-module, where $\xi$ is a partition with a hook shape. We determine the graded character of these modules. In section 5, we prove the existence of a Demazure flag for modules arising from tensor product of local Weyl module with Demazure modules of level 1 and determine the graded character formula for the tensor product of two local Weyl modules in terms of level 2 Demazure modules.
	
	\section{PRELIMINARIES}
	In this section we set the notation for the paper and recall the definitions that are willl be used throughout.
	
	\subsection{} Throughout this paper, $\mathbb{C}$ will denote the field of complex numbers, $\mathbb{Z}$ (resp. $\mathbb{Z_{+}}$), the set of integers (resp. non-negative integers) and  $\mathbb{C}[t]$, the polynomial algebra in an indeterminate $t$. For $n,r \in \mathbb{Z}_{+}$, set 
	$$[n]_{q} = \frac{1-q^n}{1-q}, \quad \begin{bmatrix} n \\ r \end{bmatrix}_{q}  = \frac{[n]_{q} [n-1]_{q} \dots [n-r+1]_{q}}{[r]_{q} [r-1]_{q} \dots [1]_{q}}$$
	$$\begin{bmatrix} n \\ 0 \end{bmatrix}_{q} =1, \quad \begin{bmatrix} n \\ 1 \end{bmatrix}_{q} = \frac{1-q^n}{1-q}, \quad \begin{bmatrix}n \\ r \end{bmatrix}_{q} = 0, \quad \text{ unless } \{n,r,n-r\}\subset \mathbb{Z}_{+}.$$

	\subsection{} Given any complex Lie algebra $\mathfrak{a}$, denote by $U({\mathfrak{a}})$ the corresponding universal enveloping algebra and by $\mathfrak{a}[t]$
	the associated current Lie-algebra with underlying vector space $\mathfrak{a}\otimes \mathbb{C}[t]$ and Lie-bracket :
	$$[x\otimes f , y\otimes g]= [x , y] \otimes fg \quad \forall x,y \in \mathfrak{a} \text{ and }f,g \in \mathbb{C}[t]$$
	The natural grading on $\mathbb{C}[t]$ defines a $\mathbb{Z}_{+}$-grading on $\mathfrak{a}[t]$ and $U(\mathfrak{a}[t])$. With respect to it, any element of $U(\mathfrak{a}[t])$ is of the form $(x_{1} \otimes t^{a_{1}})(x_{2} \otimes t^{a_{2}}) \dots (x_{n} \otimes t^{a_{n}})$ and has grade $a_{1}+a_{2}+ \dots + a_{n}$. 
	
	With respect to an ordered basis of $\mathfrak{a}[t]$, the PBW filtration on $U(\mathfrak{a}[t])$ is given as follows:
	$$U(\mathfrak{a}[t])^{\leq r} = \{a_1^{s_1}a_2^{s_2}\cdots a_l^{s_l}: 0\leq \sum\limits_{j=1}^l j s_j \leq r\}$$ We say length of $X\in U(\mathfrak{a}[t])$ is $r$ if 
	$X\in U(\mathfrak{a}[t])^{\leq r}\backslash U(\mathfrak{a}[t])^{\leq r-1}$.


	\subsection{} Let $\mathfrak{sl}_2$ be the 
	Lie algebra of two by two  trace zero matrices with entries in $\mathbb C$ and $\mathfrak{sl_2}[t]$ the associated current algebra. Let $\{x,h,y\}$ be the standard basis, $\alpha$ (resp. $\omega$) the unique simple root (resp. fundamental weight) of $\mathfrak{sl}_2$.
	
	\subsection{} For  $n\in \mathbb Z_+$, 
	let $V(n)$ be the $n+1$-dimensional irreducible $\mathfrak{sl}_2$-module generated by a highest weight vector $v_{n}$ with the following defining relations:
	$$ xv_{n} = 0, \quad hv_{n} = nv_{n}, \quad y^{n+1}v_{n} = 0.$$ 
	Any finite-dimensional $\mathfrak{sl}_2$-module $V$ can be decomposed as:
	$$V = \bigoplus_{k \in \mathbb Z} V_{k}, \qquad V_{k}= \{v \in V\mid hv=kv\}.$$ 
	The formal character of a finite-dimensional $\mathfrak{sl}_2$-module $V$, denoted as $ch_{\mathfrak {sl}_2} V$ is an element of $\mathbb{Z}[e(\omega)]$ such that 
	$$ V= \bigoplus_{k\in \mathbb Z_+} V(k)^{\oplus m_k} \qquad \Longrightarrow \qquad
	ch_{\mathfrak{sl}_2}\ V = \sum\limits_{k\in \mathbb Z_+} m_k \ ch_{\mathfrak{sl}_2} \ V(k),$$
	where $ch_{\mathfrak{sl}_2}\ V(k) = \sum\limits_{i=0}^k e((k-2i)\omega)$ for $k\in \mathbb Z_+.$
	
	\subsection{} \label{graded modules} A graded $\mathfrak{sl}_2[t]$-module $V$  is a $\mathbb{Z}$-graded vector space $V= \bigoplus\limits_{k \in \mathbb{Z}} V[k]$, which admits a $\mathfrak{sl}_2[t]$ action in the following way:
	$$(X \otimes t^r)V[k] \subseteq V[k+r], \quad r \in \mathbb{Z}_{+}, k \in \mathbb{Z}.$$
	For each $k$,  $V[k]$ is a $\mathfrak{sl}_2$-submodule of $V$. 
	If $V$ and $W$ are graded $\mathfrak{sl}_2[t]$-modules, then so is $V \otimes W$, with the following $\mathfrak{sl}-2[t]$-action:
	$$(X\otimes t^r)(v\otimes w) = (X \otimes t^r)v \otimes w + v \otimes (X \otimes t^r)w$$
	Thus in this case we have, $(V \otimes W)[k] = \bigoplus\limits_{j \in \mathbb{Z}_{+}} V[j] \otimes W[k-j]$, where $W[j] = 0 $ if $j<0$.\\
	
	\noindent Given a $\mathfrak{sl}_2$-module $W$, we define a $\mathfrak{sl}_2[t]$-structure on $W$ as follows:
	$$(X\otimes t^r)w = \delta_{r,0}Xw, $$ 
	We denote such a module by $ev_{0}(W)$. Clearly it is $\mathbb{Z}$-graded $\mathfrak{sl}_2[t]$-module with $ev_{0}(W)[k] = \delta_{k,0} W $.  For $r\in Z$, let 
	$\tau_{r}$ be the shift operater given $ (\tau_r W )[k] = W [k-r].$  Given $(n,r)\in \mathbb{Z}_{+}^2$, we denote the irreducible $\mathfrak {sl}_2[t]$-module $\tau_{r}\ ev_{0}V(n)$ by $V(n,r)$ and in particular denote  $ev_0 V(n)$ by $V(n,0)$. It was shown in \cite{MR2351379} that any irreducible, graded, finite-dimensional module is isomorphic to $V(n,r)$ for a pair $(n,r)\in \mathbb Z_+^2.$ \\
	
	\noindent The graded character of a graded $\mathfrak{sl}_2[t]$-module $V = \bigoplus\limits_{k \in \mathbb{Z}_{\geq 0}} V[k]$ is given by:
	$$ch_{gr}\ V = \sum\limits_{k\geq 0} ch_{\mathfrak{sl}_2}\ V[k]\ q^k \in \mathbb{Z}[e(\omega)][q],\qquad \text{  where } q \text{ is an indeterminate}.$$

	\subsection{} For $\lambda \in \mathbb Z_+$, the local Weyl module $W_{loc}(\lambda)$ is a finite-dimensional, $\mathbb{Z}_{+}$-graded $\mathfrak{sl}_2[t]$-module generated by an element $w_{\lambda}$ with the following defining relations:
	$$(x \otimes t^r) w_{\lambda} = 0, \quad (h \otimes t^r) w_{\lambda} = \lambda\delta_{r,0} w_{\lambda}, \quad(y \otimes 1)^{\lambda+1} w_{\lambda} = 0, \qquad \quad
	\forall r \in \mathbb{Z}_{\geq 0}.$$ By definition,  
	$W_{loc}(\lambda)$ is a highest weight representation of $\mathfrak{sl}_2[t]$ with highest weight $\lambda\omega$ and any  finite-dimensional integrable $\mathfrak{sl}_2[t]$-module with highest weight $\lambda\omega$ is a quotient of $W_{loc}(\lambda)$. In particular, $ev_{0}V(\lambda)$ is the unique irreducible quotient of $W_{loc}(\lambda)$. It was proven in \cite{MR3210603} that the graded character of $W_{loc}(m)$ is given as follows:
	
	$$ch_{gr} \ W_{loc}(m) = \sum\limits_{l=0}^{m}(\begin{bmatrix} m \\ l\end{bmatrix}_{q} - \begin{bmatrix} m \\ l-1 \end{bmatrix}_{q})\ ch_{gr}\ ev_{0}V(m-2l ).$$

	\section{PRESENTATION OF $W_{loc}(m) \otimes V(\bold{n})$}
	In \cite{MR3296163}, a family of graded representations of the current Lie algebras, which we shall refer to as CV-modules,  was introduced.
	It was shown that local Weyl modules, Demazure modules of various levels, truncated Weyl modules and fusion product modules for $\mathfrak{sl}_2[t]$ can be realized as CV-modules \cite{MR3296163, MR3407180}. In this section we recall the definition and properties of CV-modules and using them give a presentation of the $\mathfrak{sl}_2[t]$-modules that are obtained by taking tensor products of local Weyl modules with certain Demazure modules. 
	
	\subsection{} We first recall the definition of a CV-module $V(\xi)$ associated to a partition $\xi$.
	\begin{definition} Given a partition $\bon=(n_1\geq n_2\geq \cdots \geq n_k >0)$ of $n\in \mathbb Z_+$, let  $V(\bon)$ be a cyclic $\mathfrak{sl}_{2}[t]$-module generated by an element $v_{\bon}$ which satisfies the following relations:
		\begin{equation} \label{local relations}
			(x \otimes t^r) v_{\bon} = 0, \quad (h\otimes t^r) v_{\bon} = n\delta_{r,0} v_{\bon}, \quad(y \otimes 1)^{n+1} v_{\bon} = 0 
		\end{equation}  
		\begin{equation} \label{CV relations}
			(x \otimes t)^s (y \otimes 1)^{r+s} v_{\bon} = 0 \quad \forall\, r+s \geq 1+rl+\sum_{j \geq l+1} n_{j},  \text{ for some } l \in \mathbb{N} 
		\end{equation}    
		Here, $ r \in \mathbb{Z}_{+}$ in \eqref{local relations} and $r,s \in \mathbb{N}$ in \eqref{CV relations}. \end{definition} 
	Defining the grade of $v_{\bon}$ to be zero, we see that  $V(\bon)$ is a finite-dimensional graded $\mathfrak{sl}_{2}[t]$-quotient of $W_{loc}(n)$.
	The CV-modules are defined in \cite{MR3296163} for an arbitrary current algebra, but for our purposes, we need the definition only for $\mathfrak{sl}_2[t]$.\\
	
	\noindent Certain special cases of $V(\bon)$ discussed in \cite{MR3296163} are as follows:
	\begin{equation}\begin{array}{ll} \label{special cases}
			V(\bon) \cong W_{loc}(n), & \text{for }\, \bon = (1^n) , \\
			V(\bon) \cong ev_{0}(V(n)) & \text{for }\,  \bon = (n) , \\
			V(\bon) \cong D(l, lr+n_0) & \text{for }\, \bon=(l^r,n_0) ,
	\end{array}\end{equation}
	where $D(l, lr+n_0)$ denotes the Demazure module of level $l$ with highest weight $(lr+n_0)\omega = n\omega$.
	
	\subsection{} Define a power series in the indeterminate $u$ as follows: 
	\begin{eqnarray*}
		H(u) = \exp \left(- \sum_{r=1}^\infty \frac{h\otimes t^r}{r} u^r\right)
	\end{eqnarray*}
	\\
	Let $P(u)_k$ denote the coefficient of $u^k$ in  $H(u)$. 
	For $s,r \in \mathbb{Z}_{+}$, $x\in \mathfrak{g}$, define elements $x(r,s) \in U(\mathfrak{g}[t])$ as follows: 
	$$ x(r,s) = \sum\limits_{(b_{p})_{p \geq 0}} (x \otimes 1)^{(b_{0})}(x \otimes t^{1})^{(b_{1})}   \cdots (x \otimes t^{s})^{(b_{s})}$$
	with $p \in \mathbb{Z}, b_{p} \in \mathbb{Z}_{+}$ such that $ \sum\limits_{p\geq 0}b_{p} = r $, $\sum\limits_{p\geq 1}pb_{p} = s$ and for any $X \in \mathfrak{g}$, $X^{(p)} = \dfrac{X^{p}}{p!}$. 
	The following result was proved in \cite{MR0502647} and reformulated in its present form in \cite{MR1850556}.
	
	\begin{lemma}
		Given $s\in N, r \in \mathbb{Z}_{+}$, we have 
		
		\begin{equation}\label{G.eq}
			(x\otimes t)^{(s)}(y\otimes 1)^{(r+s)} - (-1)^s \big(\ \underset{k\geq 0}{\sum}\  y(r,s-k){P(u)}_k\big) \in U(\mathfrak g[t])\mathfrak n^+[t].
		\end{equation} 
		\label{G}
	\end{lemma}   
	Using \lemref{G} and \eqref{local relations}, 
	relation \eqref{CV relations} in definition of $V(\bon)$ can be replaced by the following:
	$$y(r,s) v_{\bon} =0, \text{ if } r+s \geq 1+rl+\sum_{j\geq l+1}n_{j}.$$ 
	
	\subsection{} Using the presenation of $W_{loc}(n)$ as a $CV$-module, we now give an alternate presentation of $W_{loc}(n)$ as a lowest weight module. In what follows, we use the notation : $a\otimes t^r :=a_r,$ for all $a\in \mathfrak{sl}_2, \ r\in \mathbb Z_+$.
	
	\begin{lemma} \label{W loc-} Given $n\in \mathbb Z_+$, let $W(-n)$ be the $\mathfrak{sl}_2[t]$-module generated by a vector $w_{-n}$ which satisfies the following relations:
		$$\begin{array}{lll}
			y_r w_{-n} &= 0,  \quad &\forall\, r\geq 0,\\
			h_r w_{-n} &= \delta_{0,r} (-n) w_{-n}, \quad &\forall\, r\geq 0,\\
			x(r,s) w_{-n} &=0,  \quad &\text{if } r+s \geq 1+rl+n-l .   \end{array}$$ 
		Then, $W(-n)$ is isomorphic to $W_{loc}(n).$
	\end{lemma}
	\proof Let  $W_{loc}(n)$ be generated by a vector $w_n$ which satisfies the relations \eqref{local relations} and \eqref{CV relations} for $\bon=(1^n).$ Then $y_0^nw_n$ is a non-zero vector in $W_{loc}(n)$ such that 
	$$\begin{array}{llll}
		y_ry_0^nw_n=0,\, & \forall\, r\geq 0,\quad \quad
		y(n,r).w_n = 0, &\forall\, r\geq 1.\end{array}$$
	Further, as  
	$$\begin{array}{c}
		y(n,1)= A h_1y^{n}_{0}w_n, \\
		y(n,r)w_n = A_{r}h_ry_{0}^{n}w_n +A_{r-1} h_{r-1}y(n,1)w_n
		+\cdots + {A_1 h_1 y(n,r-1)w_n} 
		, \, \forall\,  r\geq 2 \end{array} $$
	with $A, A_1,\cdots,A_r\in \mathbb Q$,
	using induction and the above relations we see that $$h_ry^{n}_{0}w_n =\delta_{r,0}(-n)w_n.$$
	Now by applying to \eqref{G.eq} the $\mathfrak{sl}_2[t]$-isomorphism that maps $y_s \rightarrow x_s$, 
	$x_s \rightarrow y_s$ and $h_s \rightarrow -h_s$
	for $s\geq 0$,  using the above relations  we get,
	\begin{equation}\label{G.eq'}
		(y\otimes t)^{(s)}(x\otimes 1)^{(r+s)}y_0^nw_{\bon} = (-1)^s   x(r,s)(-n)y_0^nw_{\bon} \end{equation} As, $x_0^{(r+s)}.y_0^nv_{\bon}=0,$ for $r+s\geq n+1$, it follows that, 
	$$  x(r,s) y_0^n w_{\bon}=0,\, \qquad  \text{ when }\hspace{.25cm} r+s\geq 1+rk+n-k.$$
	Thus, the vector $y_0^nw_n\in W_{loc}(n)$ satisfies the relations 
	$$\begin{array}{lll}
		(y \otimes t^r ) y_{0}^{n}w_{n} &= 0,  \quad &\forall\, r\geq 0,\\
		(h \otimes t^r ) y_{0}^{n}w_{n} &= \delta_{0,r} (-n) y_{0}^{n}w_{n}, \quad &\forall\, r\geq 0,\\
		x(r,s) y_{0}^{n}w_{n} &=0,  \quad &\text{if } r+s \geq 1+rl+n-l.   \end{array}$$ 
	Therefore by PBW theorem, $W_{loc}(n) = \bu (\mathfrak n^+[t])y_0^nw_n$. 
	
	Let $\phi: W(-n)\rightarrow W_{loc}(n)$ be a map such that $\phi(w_{-n})= y_0^nw_n$. Then it is clear from above that $\phi$ defines a surjective $\mathfrak{sl}_2[t]$-homomorphism from $W(-n)$ to $W_{loc}(n).$ Since by definition, $W(-n) = U(\mathfrak n^+[t])w_{-n}$, it follows from above that $\phi$ is in fact an isomorphism. \endproof
	
	\subsection{} The following results on the  $\mathfrak{sl}_{2}[t]$-modules $W{loc}(m)$ and $V(\bon)$ were proven in \cite{MR1850556} and \cite{MR3296163}. We record them for future reference. 
	\begin{lemma}\label{CV iso} Let $m,n\in \mathbb N$.
		\begin{itemize}
			\item[i.]For $m\in \mathbb N$, $\dim W_{loc}(m\omega)= 2^m.$
			\item[ii.] Given a partition $\bon=(n_1\geq n_2\geq \cdots\geq n_k)$ of $n$, let $$\bon^{+} = (n_1,n_2\cdots,n_{k-1}+1,n_k-1); \quad \bon^{-} = (n_{1} \geq n_{2}\geq \cdots \geq n_{k-1}-n_k). $$
			For $k>1$, there exists a short exact sequence of $\mathfrak{sl}_{2}[t]$-modules,
			$$0 \rightarrow \tau_{(k-1)n_k}V(\bon^{-})\xrightarrow{\phi_{-}} V(\bon)\xrightarrow{{_+}\phi} V(\bon^{+})\rightarrow 0 .$$ Furthermore, $\dim V(\bon) = \prod\limits_{i=1}^k (n_i+1)$ and $$\mathbb B(\bon)= \{v_\bon, (y_0)^{i_1 } (y_1)^{i_2 } \cdots (y_{k-1})^{i_k}v_{\bon} : (i_1 , i_2 , \dots , i_k )\in J(\bon)\}$$ is a basis of $V(\bon)$, where $$J(\bon) =\{(i_1 , i_2 , \dots , i_k)\in \mathbb Z_{\geq 0}^k: (ji_{r-1}+(j+1)i_{r})+2\sum_{l=r+1}^{n}i_{l} \leq \sum_{p=r-j}^{k} n_{p}, 
			2 \leq r \leq k+1, 1 \leq j \leq r-1 \}.$$
		\end{itemize}
	\end{lemma}
	
	\vspace{.35cm}
	
	\begin{lemma} \label{gen.W_loc.x.V(p)} Given $m, n\in \mathbb N$, let $\bon=(n_1\geq \cdots\geq n_k)$ be a partition of $n$. Let $W_{loc}(m)$ and $V(\bon)$ be the graded $\mathfrak{sl}_2[t]$-modules with generators $w_m$ and $v_\bon$, as defined above.
		Then $\mathfrak{sl}_2[t]$-module
		$W_{loc}(m\omega)\otimes V(\bon)$ is a cyclic module generated by $y_0^mw_m\otimes v_\bon$ and \begin{equation}\label{dim.W.m.n}\dim W(m\omega)\otimes V(\bon) = 2^{m}\prod\limits_{i=1}^k (n_i+1).\end{equation}.
	\end{lemma}
	\proof Given non-zero vectors $v_\bon$ and $w_m$ satisfying the given conditions, let $S$ be the submodule of $W_{loc}(m)\otimes V(\bon)$ generated by the vector 
	$y_0^mw_m\otimes v_\bon.$ Since $\mathfrak n^-.y_0^mw_m =0$ and $V(\bon)=U(\mathfrak n^-[t])v_\bon$, $$S \supseteq U(\mathfrak n^-[t])y_0^mw_m\otimes v_\bon = y_0^mw_m\otimes U(\mathfrak n^-[t])v_\bon = y_0^mw_m\otimes V(\bon). $$ 
	On the other hand, by \lemref{W loc-}, $W_{loc}(m) = U(\mathfrak n^+[t])y_0^mw_m$. Hence, for $r\in \mathbb N$ and $v\in V(\bon)$,  
	$$x_r.y_0^mw_m\otimes V(\bon)= x_r y_0^mw_m \otimes v \mod y_0^mw_m\otimes V(\bon), $$ implying $x_ry_0^mw_m\otimes V(\bon)\subseteq S$.  Now using an obvious induction argument on the length of $X\in U(\mathfrak n^+[t])$ for all basis elements $Xy_0^mw_m$ of $W_{loc}(m)$,  we see that 
	$W_{loc}(m)\otimes V(\bon)\subseteq S$. This shows that 
	$W_{loc}(m)\otimes V(\bon)$ is generated by the element $y_{0}^{m}w_{m} \otimes v_{\bon}$. 
	
	\noindent Using \lemref{CV iso}, it follows from \ref{graded modules} that $W_{loc}(m)\otimes V(\bon)$ is a finite-dimensional, graded $\mathfrak{sl}_{2}[t]$-module such that \eqref{dim.W.m.n} holds. \endproof

	\subsection{} Let $\mathfrak b^\pm$ be the subalgebra  $\mathfrak h+\mathfrak n^\pm$ of $\mathfrak g$. Then $\mathfrak g= \mathfrak b^+\oplus\mathfrak n^-
	=\mathfrak b^-\oplus\mathfrak n^+$ and by the Poincar\'e-Birkhoff-Witt theorem, 
	\begin{equation} \label{U decompose}
		U(\mathfrak g[t])=U(\mathfrak n^-[t])\oplus U(\mathfrak g[t])\mathfrak b^+[t]=U(\mathfrak n^+[t])\oplus U(\mathfrak g[t])\mathfrak b^-[t].\end{equation}
	Let $\operatorname{\bold{pr}}^\pm:U(\mathfrak g[t])\rightarrow U(\mathfrak n^\pm[t])$ be the projections with respect to decompositions \eqref{U decompose}. 
	
	\begin{lemma}\label{Ideals} Let 
		$\bon=(n_1\geq n_2\geq \cdots\geq n_k)$ be a partition of a natural number $n$ and let $V(\bon)$ be the associated $CV$-module generated by a highest weight vector $v_\bon$. Let $J_{\bon}$ be the ideal of $U(\mathfrak n^-[t])$ generated by the elements 
		\begin{equation} \label{gen J_n}\{\pr^-((x_1)^s(y_0)^{r+s}): s+r\geq 1+rl+\sum\limits_{j\geq l+1} n_j,\, \text{ for some}\ k\in \mathbb N\}.\end{equation} 
		The map $U(\mathfrak n^-[t])\rightarrow V(\bon)$ which sends $g$ to $gv_\bon$ induces a vector space isomorphism \begin{equation} \label{U J.bon ideal} U(\mathfrak n^-[t])/J_\bon \cong V(\bon).\end{equation}
	\end{lemma}
	\proof Given the partition $\bon$ of $n$, for all  $Y\in U(\mathfrak n^-[t]\oplus \mathfrak b^+t[t])$ we have, $$Yv_\bon=\pr^-(Y)v_\bon.$$  Thus by definition of the CV-module $V(\bon)$, the kernel of the map $U(\mathfrak n^-[t])\rightarrow V(\bon)$ which sends $g$ to $gv_\bon$  is generated by the set \eqref{gen J_n}.
	
	\noindent Using the defining relations of the $CV$-module $V(\bon)$ we 
	have $$y(r,s)v_\bon =0, \quad \text{ for } r+s\geq 1+r(k-1) + n_k.$$ 
	Hence $y_{k-1}^{i_k}v_\bon =0$ whenever $i_k> n_k$. 
	
	\noindent Furthermore, as $\bon^+$ is the partition $(n_1^+\geq \cdots\geq n_{k-1}^+\geq n_{k}-1)$ with $$n_{j_0}^+ = n_{k-1}+1, \quad \text{ if }  j_0 = \min\{ j : n_j=n_k-1 , 1\leq j\leq k-1\}$$ and $n_j^+ = n_j$ for $1\leq j\leq k-1$, $j\neq j_0$, 
	using the short exact sequence in \lemref{CV iso}(b) repeatedly, we see that for $0\leq i_k<n_k$, if $\bon^+(i_k)=(\bar{n}_1\geq \cdots\geq \bar{n}_{k-1}\geq i_k)$ is a partition of $n$ such that ${{_+}\phi}^{n_k-i_k}(V(\bon)) = V(\bon^+(i_k))$, then from the definition of the map ${_+}\phi$, it follows that 
	$$\bar{n}_{k-1} \leq \min \{n_{k-1}+n_k-i_{k}, n_j+\frac{n_k-i_k -\sum\limits_{s=1}^{k-j}(s)(n_{k-s-1}-n_{k-s})}{j}, 1\leq j\leq k-2\}.$$ 
	Therefore, using the short exact sequence for $V(\bon^+(i_k))$, we see that 
	$$\begin{array}{rl}y_{k-2}^{i_{k-1}}y_{k-1}^{i_k}v_\bon =0, \quad& \text{ if } i_{k-1}> \bar{n}_{k-1}-i_k,\\ i.e., \qquad \quad y_{k-2}^{i_{k-1}}y_{k-1}^{i_k}v_\bon =0, \quad& \text{ if } ji_{k-1}+(j+1)i_k>\sum\limits_{r=k-j}^k n_r, \quad \text{ for some } j\leq k-1.\end{array}$$ 
	As $\sum\limits_{j=1}^{k-1}\bar{n}_j=\sum\limits_{j=1}^k n_j-i_k$, using the arguments recursively, we see that for a $k$-tuple $(s_1,\cdots,s_k)$, 
	$$y_{r-1}^{s_{r-1}}y_r^{s_r}\cdots y_{k-1}^{s_k}v_\bon =0 \quad \text{ if } (ji_{r-1}+(j+1)i_{r})+2\sum\limits_{l=r+1}^{n}i_{l} > \sum\limits_{p=r-j}^{k} n_{p},$$ for some $ 0<r-1\leq k-1$ and $1\leq j\leq r-1$. This shows that all monomials of the form $y_0^{s_1}\cdots y_{r-1}^{s_{r-1}}y_r^{s_r}\cdots y_{k-1}^{s_k}\in J_\bon$, whenever  $(ji_{r-1}+(j+1)i_{r})+2\sum\limits_{l=r+1}^{n}i_{l} > \sum\limits_{p=r-j}^{k} n_{p}$ for some $0<r-1\leq k-1$ and $1\leq j\leq r-1$. Consequently, $U(\mathfrak n^-[t])/J_\bon$ is spanned by cosets of the elements $(y_0)^{i_1 } (y_1)^{i_2 } \cdots (y_{k-1})^{i_k}$ with $(i_1, i_2, \cdots, i_k)\in J(\bon)\cup (0,\cdots,0).$ 
	
	On the other hand by \lemref{CV iso}, $\{v_\bon, (y_0)^{i_1 } (y_1)^{i_2 } \cdots (y_{k-1})^{i_k}v_{\bon} : (i_1 , i_2 , \dots , i_k )\in J(\bon)\}$ is a basis of $V(\bon)$. Thus it follows that the map $U(\mathfrak n^-[t])\rightarrow V(\bon)$ gives a vector space isomorphism between $U(\mathfrak n^-[t])/J_\bon$ and $V(\bon)$ . \endproof
	
	\begin{remark}\label{Ideal-Im} Similarly,  it can be shown that if $I_m$ is the ideal of $U(\mathfrak n^+[t])$ generated by the set 
		$\{\pr^+((y_1)^k(x_0)^\ell): \ell > m\}$, 
		then the map $U(\mathfrak n^+[t])\rightarrow W_{loc}(m)$ given by $g\mapsto gy_0^mw_m$ gives a vector space isomorphism between  $U(\mathfrak n^+[t])/I_m$ and $W_{loc}(m)$. In this case $U(\mathfrak n^+[t])/I_m$ is spanned by the set 
		$$\{(x_0)^{i_1 } (x_1)^{i_2 } \cdots (x_{m-1})^{i_m} : (i_1, i_2, \cdots, i_m)\in J(1^n)\cup (0,\cdots,0)\}.$$ \end{remark}
	
	\vspace{.25cm}
	
	\subsection{The module {$W(m,\bon)$}} Given $n,m\in \mathbb N$ and a partition  $\bon=(n_1\geq n_2\geq \cdots \geq n_k)$ of $n$, define $W(m,\bon)$ to be the $\mathfrak{sl}_{2}[t]$-module generated by an element $w_{\bon}^{m}$ which satisfies the following relations:
	\begin{align} 
		(h \otimes t^r )w_{\bon}^{m} &= (n-m) \delta_{r,0} w_{\bon}^{m}\label{h reln}\\
		x(r,s)w_{\bon}^{m} &= 0,\text{ if } r+s \geq 1+rl+m-l \text{ for some }l \in \mathbb{N}\label{x reln}\\
		y(r,s)w_{\bon}^{m} &= 0, \text{ if } r+s \geq 1+rl+\sum_{j \geq l+1} n_{j} \text{ for some }l \in \mathbb{N} \label{y reln}
	\end{align}
	
	\noindent The following is the main result of this section : 
	\begin{theorem} \label{present} Given $n,m\in \mathbb N$ and $\bon=(n_1\geq n_2\geq \cdots \geq n_{k} )$, a partition of $n$, the  $\mathfrak{sl}_2[t]$-module $W(m,\bon)$ is isomorphic to $W_{loc}(m)\otimes V(\bon)$.
	\end{theorem}
	
	\noindent Observe that \thmref{present} gives a presentation of the $\mathfrak{sl}_2[t]$-module $W_{loc}(m)\otimes V(\bon)$. The following lemmas help prove the theorem. 
	
	\begin{lemma} \label{lemma1}
		There exists a surjective, $\mathfrak{sl}_{2}[t]$-module homomorphism from $W(m,\bon)$ to $W_{loc}(m)\otimes V(\bon)$.
	\end{lemma} 
	\proof By \lemref{gen.W_loc.x.V(p)}, we know that 
	if $w_m$ and $v_\bon$ are the highest weight generators of $W_{loc}(m)$ and $V(\bon)$ respectively, then $y_0^mw_m\otimes v_\bon$ generates $W_{loc}(m)\otimes V(\bon).$
	
	Define a map $\phi: W(m,\bon) \rightarrow W_{loc}(m)\otimes V(\bon)$, such that, $\phi(w_{\bon}^{m})= y_{0}^{m}w_{m} \otimes v_{\bon}$. Since  $y_0^mw_m\otimes v_\bon$ generates $W_{loc}(m)\otimes V(\bon)$ and for $r\geq 0$, using \lemref{W loc-} we have, 
	\begin{align*}
		(h \otimes 1)(y_{0}^{m}w_{m} \otimes v_{\bon}) &= h_0 y_{0}^{m}w_{m} \otimes v_{\bon}+ y_{0}^{m}w_{m} \otimes h_0 v_{\bon} = (n-m)(y_{0}^{m}w_{m} \otimes v_{\bon})\\
		x(r,s)(y_{0}^{m}w_{m} \otimes v_{\bon})&= x(r,s)y_{0}^{m}w_{m} \otimes v_{\bon} = 0  \quad \text{ if } r+s \geq 1+rl + m-l \\
		y(r,s)(y_{0}^{m}w_{m} \otimes v_{\bon})&= y_{0}^{m}w_{m} \otimes y(r,s)v_{\bon} = 0 \quad \text{ if } r+s \geq 1+rl + \sum_{j \geq l+1} n_{j} \quad (\text{By \eqref{CV relations}}) 
	\end{align*}  and for $r\geq 1$,
	\begin{align*}
		(h \otimes t^r )(y_{0}^{m}w_{m} \otimes v_{\bon}) &= h_{r}y_{0}^{m}w_{m} \otimes v_{\bon} = 0. \quad 
	\end{align*} 
	This shows that $\phi$ is a well-defined surjective $\mathfrak {sl}_2[t]$-module homomorphism. \endproof

	\begin{lemma} \label{lemma2}
		$dim(W(m,\bon)) \leq 2^{m}\prod\limits_{j=1}^k(n_j+1)$.
	\end{lemma}
	\proof By definition, $W(m,\bon) = U(\mathfrak{sl}_{2}[t])w_{\bon}^{m}$, where $w_\bon^m$ is a non-zero vector satisfying the relations \eqref{h reln}-\eqref{y reln}. Using Poincar\'e-Birkhoff-Witt theorem and \eqref{h reln}, we can write $$U(\mathfrak{sl}_{2}[t]) w_{\bon}^{m} = U(\eta^{+}[t])U(\eta^{-}[t])U(\mathfrak{h}[t])w_{\bon}^{m} = U(\eta^{+}[t])U(\eta^{-}[t])w_{\bon}^{m}. $$
	For the partition $\bon$, let $J_\bon$ be the ideal of $U(\mathfrak n^-[t])$ as 
	described in \lemref{Ideals}. As a consequence of the defining relations 
	\eqref{h reln} and \eqref{y reln}, we see that the same arguments as those 
	used in \lemref{Ideals} show that $Yw^m_\bon=0$ for all $Y\in J_\bon$. Hence it 
	follows that if 
	$${\mathcal B}_\bon = \{y_{0}^{i_{1}}y_{1}^{i_{2}} \dots y_{k-1}^{i_{k}}| (i_{1}, i_{2}, \dots , i_{k}) \in J(\bon)\cup (0,0,\cdots,0)\},$$ then $$W(m,\bon) \subseteq \sum_{Y\in {\mathcal B}_\bon } U(\mathfrak n^+[t])Yw^m_\bon. $$
	
	\noindent For $Y\in U(\mathfrak n^-[t]$, let $\ell{n}(Y)$ be the PBW length of $Y$. \\
	For $Y\in U(\mathfrak n^-[t])$, when $\ell{n}(Y) =0$, we have,
	$U(\mathfrak n^+[t])Yw^m_\bon= U(\mathfrak n^+[t])w^m_\bon.$ In this case, as a consequence of the defining relations \eqref{h reln} and \eqref{x reln}, the similar arguments to those used in \lemref{Ideals} show that if $I_m$ is the ideal of $U(\mathfrak n^+[t])$ defined in \remref{Ideal-Im}, then, 
	as vector spaces $U(\mathfrak n^+[t])/I_m$ is isomorphic to $U(\mathfrak n^+[t])w^m_n$. Hence $\dim U(\mathfrak n^+[t])w^m_\bon\leq 2^m.$ 
	
	Given $Y\in {\mathcal B}_\bon$, by applying induction on $\ell{n}(Y)$, we now prove that
	\begin{equation}\label{ell{n}(Y)} 
		U(\mathfrak n^+[t])Yw^m_\bon \subseteq YU(\mathfrak n^+[t])w^m_\bon + \sum\limits_{\{Y'\in {\mathcal B}_\bon: \ell{n}(Y')<\ell{n}Y\}} Y' U(\mathfrak n^+[t])w^m_\bon.\end{equation}
	When $\ell{n}(Y)=1$, $Y=y_j$ with $0\leq j\leq m-1.$ In this case, given $x^{a_0}_{0}x^{a_1}_{1} \dots x^{a_p}_{p}\in U(\mathfrak n^+[t])$, we have, 
	$$
	x^{a_0}_{0}x^{a_1}_{1} \dots x^{a_p}_{p}y_j w_{\bon}^{m} = x^{a_0}_{0}x^{a_1}_{1} \dots x^{a_{p}-1}_{p}y_j x_{p} w_{\bon}^{m} + x^{a_0}_{0}x^{a_1}_{1} \dots x^{a_p -1}_{p}h_{j+p} w_{\bon}^{m}.$$
	Since $x_rw_\bon^m =0$ for $r\geq m$, we see that $x^{a_0}_{0}x^{a_1}_{1} \dots x^{a_p}_{p}y_j x_{p} w_{\bon}^{m}=0$ whenever $p\geq m$. Therefore, assuming $p\leq m-1$, we get,
	\begin{align*}
		x^{a_0}_{0}x^{a_1}_{1} \dots x^{a_p}_{p}y_j w_{\bon}^{m} &=  x^{a_0}_{0}x^{a_1}_{1} \dots x^{a_{p}-1}_{p}y_j x_{p} w_{\bon}^{m} + x^{a_0}_{0}x^{a_1}_{1} \dots x^{a_p -1}_{p}h_{j+p} w_{\bon}^{m},\\
		&= x^{a_0}_{0}x^{a_1}_{1} \dots x^{a_{p}-2}_{p}y_j x^{2}_{p} w_{\bon}^{m} + x^{a_0}_{0}x^{a_1}_{1} \dots x^{a_p -2}_{p}h_{j+p}x_{p} w_{\bon}^{m} + x^{a_0}_{0}x^{a_1}_{1} \dots x^{a_p -1}_{p}h_{j+p} w_{\bon}^{m}\\
		& \vdots \\
		& = y_{j}x^{a_0}_{0}x^{a_1}_{1} \dots x^{a_p}_{p} w_{\bon}^{m} + u, \quad\text{where } u \in U(\mathfrak n^{+}[t])w_{\bon}^{m}.
	\end{align*} 
	Thus \eqref{ell{n}(Y)} holds $\ell{n}(Y)=1$. Now assume that \eqref{ell{n}(Y)} holds for all $Y\in {\mathcal B}_\bon$ with  $\ell{n}(Y)<r$ and consider $ x^{a_0}_{0}x^{a_1}_{1} \dots x^{a_p}_{p}Yw_\bon^m$, with $Y=y^{i_{1}}_{0}y^{i_{2}}_{1} \dots y^{i_{k}}_{k-1}$ with $\ell{n}(Y)=r$. 
	Then,$$ \begin{array}{l}
		x^{a_0}_{0}x^{a_1}_{1} \cdots x^{a_p}_{p}y^{i_{1}}_{0}\cdots y^{i_{k}}_{k-1}w_{\bon}^{m}\\= (x^{a_0}_{0}x^{a_1}_{1} \cdots x^{a_p -1}_{p}y_{0}x_{p}y^{i_{1} -1}_{0}y^{i_{2}}_{1} \dots y^{i_{k}}_{k-1} + x^{a_0}_{0}x^{a_1}_{1} \dots x^{a_p -1}_{p}h_{p}y^{i_{1} -1}_{0}y^{i_{2}}_{1} \dots y^{i_{k}}_{k-1}) w_{\bon}^{m} \\
		=x^{a_0}_{0}x^{a_1}_{1} \cdots x^{a_p -1}_{p}y_{0}x_{p}y^{i_{1} -1}_{0}y^{i_{2}}_{1} \dots y^{i_{k}}_{k-1}w_{\bon}^{m}  + (i_1-1){ x^{a_0}_{0}x^{a_1}_{1} \dots x^{a_p -1}_{p}}y_0^{i_1-2}y^{i_{2}}_{1} \dots y_p^{i_p+1}\cdots y^{i_{k}}_{k-1}w_{\bon}^{m}\\+\sum\limits_{j=2}^{k-p-1}(i_j){ x^{a_0}_{0}x^{a_1}_{1} \cdots x^{a_p -1}_{p}} y^{i_{1}-1}_{0}y^{i_{2}}_{1} \dots y^{i_j-1}\dots y_{j+p}^{i_{j+p}+1}\cdots y^{i_{k}}_{k-1}w_{\bon}^{m}. 
	\end{array} $$ 
	Using the vector space isomorphism between $U(\mathfrak n^-[t])w_\bon^m$ and $U(\mathfrak n^-[t])/J_\bon$, it follows from \lemref{Ideals} that,
	$$x^{a_0}_{0}x^{a_1}_{1} \dots x^{a_p}_{p}y^{i_{1}}_{0}\cdots y^{i_{k}}_{k-1}w_{\bon}^{m}
	= x^{a_0}_{0}x^{a_1}_{1} \dots x^{a_p -1}_{p}y_{0}x_{p}y^{i_{1} -1}_{0}y^{i_{2}}_{1} \dots y^{i_{k}}_{k-1}w_{\bon}^{m} + u_{r-1},$$
	where $u_{r-1} \in \sum\limits_{\{ Y'\in \mathcal B_\bon: \ell{n}(Y')=r-1\} } U(\eta^{+}[t])\cdot Y' w_{\lambda}^{m}$. Thus by induction hypothesis, $$u_{r-1} \in \sum\limits_{Y\in \mathcal B_\bon } Y \cdot U(\eta^{+}[t]) \cdot w_{\bon}^{m}.$$ Now, repeating the proceedure we get, 
	$$x^{a_0}_{0}x^{a_1}_{1} \dots x^{a_p}_{p}y^{i_{1}}_{0}y^{i_{2}}_{1} \dots y^{i_{k}}_{k-1}w_{\bon}^{m}  = x^{a_0}_{0}x^{a_1}_{1} \dots x^{a_p -2}_{p}y_{0}x^{2}_{p}y^{i_{1} -1}_{0}y^{i_{2}}_{1} \dots y^{i_{k}}_{k-1}w_{\bon}^{m} + u_{r-2}+u_{r-1},$$
	where $u_{r-2}\in \sum\limits_{\{Y\in \mathcal B_\bon\}} Y U(\mathfrak n^+[t])w_\bon^m.$
	Continuing this way, we get 
	$$x^{a_0}_{0}x^{a_1}_{1} \dots x^{a_p}_{p}y^{i_{1}}_{0}y^{i_{2}}_{1} \dots y^{i_{k}}_{k-1}w_{\bon}^{m}= y_{0}x^{a_0}_{0}x^{a_1}_{1} \dots x^{a_p}_{p}y^{i_{1} -1}_{0}y^{i_{2}}_{1} \dots y^{i_{k}}_{k-1}w_{\bon}^{m} + u, $$
	where $u \in \sum\limits_{\{Y'\in \mathcal B_\bon\}}Y' \cdot U(\eta^{+}[t])w_{\bon}^{m})$. This shows that
	$W(m,\bon)$ is spanned by elements of the form $\{YXw_\bon^m: X\in \mathcal B_{{\bf{1}}^m} \,  Y\in\mathcal B_\bon\}$. Hence it follows that
	$$\qquad \qquad \qquad \dim W(m,\bon) \leq |\mathcal B_{{\bf{1}}^m}| |\mathcal B_{\bon}| = 2^m \prod\limits_{j=1}^k(n_j+1). \qquad \qquad \qquad \qquad \qquad \qed $$  
	
	\noindent We conclude this section with the proof of \thmref{present}. 
	\noindent \proof 
	From \lemref{lemma1} and \lemref{CV iso}, it follows that $\dim W(m,\bon)\geq 2^{m}\prod\limits_{i=1}^k(n_i+1). $\\ On the other hand by \lemref{lemma2} , 
	$\dim W(m,\bon) \leq 2^{m}\prod\limits_{i=1}^k(n_i+1).$\\
	Hence it follows that $\dim W(m,\bon)\geq 2^{m}\prod\limits_{i=1}^k(n_i+1)= \dim W(m)\otimes V(\bon) $.\\ This completes the proof of the theorem.  \qed

	
	\section{$CV$-filtration and Graded character of $W_{loc}(m \omega)\otimes V(n)$ }

\label{multiplicity} Given a finite-dimensional graded $\mathfrak{sl}_2[t]$-module $U$, we say $U$ admits a $CV$-filtration, if there exist an 
decreasing chain of $\mathfrak{sl}_2[t]$-submodules
$$U =U_0\supset  U_1 \supset U_2 \supset U_3 \supset \dots \supset U_r \supset (0) $$ such that for each $i$, $\frac{U_i}{U_{i+1}}\cong V(\xi^i)$, where  $\xi^i$ is a partition. On the other hand, a filtration of a representation whose successive quotients are isomorphic to Demazure modules is called an excellent filltration.  It is known that if  $V=V_1\supseteq V_2\supseteq \cdots \supseteq V_n\supset 0$ is a chain of graded $\mathfrak g[t]$-modules whose successive quotients have excellent filtration then $V$ has an excellent filtration. 
In this context the following result was proved in \cite[Theorem 3.3]{MR3210603}.
\begin{proposition}\label{Theorem 3.3} Let $n\in \mathbb N$. Given a partition $\bon=(n_1\geq n_2\geq \cdots \geq n_k)$ of $n$, the $\mathfrak {sl}_2[t]$-module $V(\bon)$ has a Demazure flag of level $l$ if and only if $n_1\leq l.$
\end{proposition}
We use this result to show that the $\mathfrak{sl}_2[t]$-modules of the form $W_{loc}(m)\otimes V(n)$ admit a Demazure flag of level $k$ for $k\geq n+1$. We 
also deduce the fermionic version of the character formula for $W_{loc}(m)\otimes V(n)$. 

\subsection{} Given $k,r\in \mathbb N$, with $k\geq r$, let $\hook(k,r):=(k-r+1,1^{r-1})$ denote the hook type partition of $k$ into $r$ parts. In this section we prove the following:
\begin{theorem} \label{filtration in irreducible}
	Let $m, n\in \mathbb N$. The $\mathfrak{sl}_{2}[t]$-module, $W_{loc}(m \omega)\otimes V(n)$ admits a $CV$-filtration whose successive quotients are isomorphic to  
	\begin{align*}
		V(\hook(m-n+2i,m-n+i))&, \quad 0 \leq i \leq n, \quad \text{when } m> n; \\
		V(\hook(n-m+2i,i), V(n-m)&,  \quad 1 \leq i \leq m, \quad \text{when } m\leq n.
	\end{align*} 
\end{theorem}

\noindent We need the following lemmas to prove the theorem. 

\begin{lemma} \label{A(r,i)} Let $r,i\in \mathbb N$.\\
	If $r \geq i $, then the $i\times i$ matrix $A(r,i)$ given by 
	$$A(r,i) = \begin{bmatrix}
		1 & 1  & \cdots  & 1\\
		{{r} \choose {1}} & {{r+1} \choose {1}}   & \cdots  & {{r+i-1} \choose {1}} \\
		\vdots & \vdots & \cdots & \vdots \\
		{{r} \choose {i-1}}  & {{r+1} \choose {i-1}}  & 
		\cdots & {{r+i-1} \choose {i-1}} 
	\end{bmatrix}$$ is invertible.\\
	If $r < i$, then the $i\times i$ matrix $B(r,i)$ given by
	$$ B(r,i) = \begin{bmatrix}
		1 & 1 & 1 & \cdots & 1& 1 & \cdots  & 1\\
		{{r} \choose {1}} & {{r+1} \choose {1}}& {{r+2} \choose {1}} & \cdots & {{i-2} \choose {1}} & {{i-1} \choose {1}}& \cdots  & {{r+i-1} \choose {1}} \\
		\vdots & \vdots & \vdots & \cdots & \vdots & \vdots & \cdots & \vdots \\
		{{r} \choose {r}} &{{r+1} \choose {r}}&{{r+2} \choose {r}} & \cdots &{{i-2} \choose {r}} & {{i-1} \choose {r}} & \cdots & {{r+i-1} \choose {r}}  \\
		0 &{{r+1} \choose {r+1}} & {{r+2} \choose {r+1}}& \cdots & {{i-2} \choose {r+1}} & {{i-1} \choose {r+1}} & \cdots & {{r+i-1} \choose {r+1}}  \\
		\vdots & \vdots & \vdots & \cdots & \vdots & \vdots & \cdots & \vdots \\
		0 & 0 & 0 & \cdots & 0 & {{i-1} \choose {i-1}}  & \cdots & {{r+i-1} \choose {i-1}} \end{bmatrix}$$ is invertible.
\end{lemma}
\proof We prove the lemma by applying induction on $i$. \\
When $i = 1$, $A(r,1) = [1]$ which is invertible.\\ 
When $i=2$,  $A(r,2) = \begin{bmatrix}
	1 & 1\\
	r & r+1
\end{bmatrix}$ and $\det A(r,2) = 1$, impying that $A(r,i)$ is invertible.\\
Assume that the result holds for all $k<i$ and consider the matrix
$A(r,i)$ 
	. On applying column transformation to $A(r,i)$ we can get the matrix 
	$$\tilde{A}(r) = \begin{bmatrix}
		1 & 0  & \cdots  & 0\\
		{{r} \choose {1}} & 1   & \cdots  & 1\\
		{{r} \choose {2}} &{{r} \choose {1}} & \cdots & {{r+i-2} \choose {1}}  \\
		\vdots & \vdots & \cdots & \vdots \\
		{{r} \choose {i-1}}  & {{r} \choose {i-2}}  & \cdots & {{r+i-2} \choose {i-2}} 
	\end{bmatrix}$$ Clearly $\det A(r,i) = \det \tilde{A}(r) =\det A(r,i-1).$ By induction hypothesis the matrix $A(r,i-1)$ is invertible. Hence, it follows that $A(r,i)$ is invertible. \\
	Similar arguments show that the matrix $B(r, i)$ is also invertible when $r < i$.
	\endproof
	
	\begin{lemma} \label{sum dim Ui/Ui-1} Let $m,n\in \mathbb N$.\\
		(i) For $m> n$,  
		$\sum\limits_{i=0}^{n} \dim V(\hook(m-n+2i, m-n+i) = 2^m(n+1)$. \\
		(ii) For $m\leq n$, 
		$\sum\limits_{i=1}^{m} \dim V(\hook(n-m+2i, i)+ \dim V(n-m) =2^m(n+1)$ . 
	\end{lemma}
	\proof  By \lemref{CV iso},
	$$\dim V(\hook(m+n-2i,m-i)) = 2^{m-i-1}(n-i+1+1) = 2^{m-i-1}(n-i+2).$$
	
	\noindent (i) When $m > n$, we have, 
	\begin{align*}
		\sum\limits_{i=0}^{n} \dim V(\hook(m-n+2i, m-n+i) & = \sum\limits_{i=0}^n (i+2)2^{m-n+i-1} = 2^{m-n-1} \sum\limits_{i=0}^n (i+2)2^i \end{align*}
	Putting $s_n = \sum\limits_{i=0}^n (i+2)2^i,$ we see that 
	$s_n=\sum\limits_{i=1}^n i2^i +\sum\limits_{i=1}^{n+1}2^i.$\\
	As $\sum\limits_{i=1}^n i2^i +(n+1)2^{n+1}-2(\sum\limits_{i=1}^n i2^{i})=2(2^{n+1}-1)$, 
	$$s_n=(n+1)2^{n+1}-2(2^{n+1}-1)+\sum\limits_{i=1}^{n+1}2^i=(n+1)2^{n+1}. $$
	Hence, \, $\sum\limits_{i=0}^{n} \dim V(\hook(m+n-2i,m-i)=2^{m-n-1}s_n = 2^m(n+1).$\\
	
	\noindent (ii) When $m\leq n$, we have,  
	\begin{align*}
		\sum\limits_{i=1}^{m} \dim V(\hook(n-m+2i,i)) + \dim V(n-m) = \sum_{i=1}^{m} (n-m+2+i)2^{i-1} + (n-m+1) \\=
		\sum_{i=0}^{m-1} 2^i (n-m+3+i) + (n-m+1)
		= (n-m+1) \sum_{i=0}^{m-1} 2^i + \sum_{i=0}^{m-1} (i+2)2^i + (n-m+1) \end{align*}
	Using part(i), we thus get,
	\begin{align*}&\sum\limits_{i=0}^{m-1} \dim V(\hook(n+m-2i,m-i)) + \dim V(n-m)\\
		&=(n-m+1)(2^m-1)+m2^m+(n-m+1)=2^m(n+1). \end{align*}
	This completes the proof of the lemma. \endproof

	\noindent We now prove \thmref{filtration in irreducible}.
	\proof Let $w_m$ and $v_n$ be the highest weight generators of the modules $W_{loc}(m)$ and $V(n)$ respectively. By \lemref{gen.W_loc.x.V(p)},  
	$W_{loc}(m)\otimes V(n)= U(\mathfrak{sl}_2[t]).w_m\otimes y^{(n)}v_n$.\\
	
	\noindent When $m> n \geq 1$,  set $$U_i=W_{n-i}= U(\mathfrak{sl}_{2}[t])(w_{m}\otimes y_{0}^{(n-i)}v_n), \qquad 0 \leq i \leq n.$$ As $x_r.w_m=0$ for all $r\in \mathbb Z_+$, using representation theory of $\mathfrak{sl}_2(\mathbb C)$ \cite[Section 7]{Humphreys}, we have 
	$$x_0.(w_m \otimes y_0^{(i)} v_n) = w_m \otimes x_0 y_0^{(i)} v_n = (n-i+1) w_m \otimes y_0^{(i-1)} v_n.$$ 
	generates $W_{loc}(m\omega) \otimes V(n)$. Hence, 
	$$W_{loc}(m\omega) \otimes V(n)=U_0 \supset U_1 \supset \cdots U_n\supset 0, $$
	is an decreasing chain of submodules of $W_{loc}(m\omega)\otimes V(n)$. To prove the result we show: 
	$$\frac{U_i}{U_{i+1}} \cong V(\hook(m-n+2i,m-n+i) = V(i+1, 1^{m-n+i-1}), \qquad \text{ for } \, 0 \leq i \leq n,$$ or equivalently,
	$\dfrac{W_i}{W_{i-1}} \cong V(\hook(m+n-2i,m-i) = V(n+1-i, 1^{m-i-1}) $ 
	for $0 \leq i \leq n.$\\

	\noindent By definition (\secref{graded modules}), 
	$$ g_r. v=0, \qquad \text{ for all } \, v\in V(n),\,  g\in \mathfrak g, \, r\in \mathbb N.$$
	Hence for $r,s\in \mathbb Z_+$, 
	\begin{equation}\label{y(r,s)(i)}
		\begin{array}{ll}
			(-1)^s y(r,s)(w_m\otimes y_0^{(i)}v_n) & = x_1^{(s)}y_0^{(r+s)}(w_m\otimes y_0^{(i)}v_n)\\ 
			& = \frac{x_1^{(s)}}{(r+s)!i !} \left(\sum\limits_{l=0}^{r+s} {r+s\choose l} y_0^{r+s-l}w_m\otimes y_0^{l+i}v_n\right)\\
			& = \, \, x_1^{(s)} \quad \left(\sum\limits_{l=0}^{r+s} {l+i\choose l} x_1^{(s)} y_0^{(r+s-l)}w_m\otimes y_0^{(l+i)}v_n\right)\\ &= \sum\limits_{l=0}^{r+s}\, \,  {l+i\choose l}\,\, x_1^{(s)}y_0^{(r+s-l)}w_m\otimes y_0^{(l+i)} v_n\\ & =
			\sum\limits_{l=0}^{r-1}\,\, {l+i\choose l}\, \,  y(r-l,s)w_m\otimes y_0^{(l+i)}v_n \end{array}\end{equation}
	Further, as $y(r,s)w_m=0$ whenever $r+s> m$,  
	$$y(r,s)(w_m\otimes y_0^{(i)}v_n) =0 \qquad \forall\, s\geq m,\, r\geq 1,\, i\geq 0 .$$
	Therefore, 
	$$ \begin{array}{c}
		x_r .(w_m \otimes v_n)= 0, \qquad  h_r.(w_m \otimes v_n) = \delta_{r,0}(n+m)(w_m \otimes v_n), \qquad \forall\, r\in \mathbb Z_+,\\
		\\
		\qquad y(r,m)(w_m \otimes v_n) = 0 \quad \forall\, 1\leq r \leq n+1,
	\end{array}$$ 
	implying that $W_0=U_n$ is a quotient of $V(\hook(m+n,m)).$
	
	\noindent Now, using standard $\mathfrak{sl}_2$-representation theory (\cite[Section 7]{Humphreys}) we see that for $i>0$,
	$$\begin{array}{ll}
		x_r.(w_m \otimes y_{0}^{(i)}v_n ) = \delta_{r,0}(n-i+1) \,  w_m \otimes y_{0}^{(i-1)}v_n \in U_{n-i+1} \\
		h_r.(w_m \otimes y_{0}^{(i)}v_n ) = \delta_{r,0}(m+n-2i) \,  w_m \otimes y_{0}^{(i)}v_n,
	\end{array}$$ 
	and putting $s=m-i$ in \eqref{y(r,s)(i)}, we get,
	\begin{equation} \label{r<i}y(r,m-i)(w_m\otimes y_0^{(i)}v_n)=\sum\limits_{l=0}^{r-1}\,\, {l+i\choose l}\, \,  y(r-l,m-i)w_m\otimes y_0^{(l+i)}v_n.\end{equation}
	Once again using the relation $y(r,s)w_m=0$ whenever $r+s>m$, we see that
	for $r\geq i$, \eqref{r<i} reduces to:
	\begin{equation} \label{r>i} \begin{array}{ll}
			(x \otimes t)^{(m-i)}(y \otimes 1)^{(r+m-i)}(w_m \otimes y_{0}^{(i)}v_n )&= \sum\limits_{l=0}^{r-1}\,\, {l+i\choose l}\, \,  y(r-l,m-i)w_m\otimes y_0^{(l+i)}v_n \\
			&=\sum\limits_{l=r-i}^{r-1}\,\, {l+i\choose l}\, \,  y(r-l,m-i)w_m\otimes y_0^{(l+i)}v_n\\
			& = \sum\limits_{l=0}^{i-1} {{r+l} \choose {i}} y(i-l, m-i)w_m \otimes y_{0}^{(r+l)}v_n );
	\end{array}\end{equation}
	and for  $r<i$, putting $r=i-p$, \eqref{r<i} reduces to:
	\begin{equation} \label{r=i-p}
		\begin{array}{ll}y(r,m-i)(w_m\otimes y_0^{(i)}v_n)&=\sum\limits_{l=0}^{r-1}\,\, {l+i\choose l}\, \,  y(i-p-l,m-i)w_m\otimes y_0^{(l+i)}v_n\\
			&=\sum\limits_{k=p}^{r+p-1}\,\, {k-p+i\choose i}\, \,  y(i-k,m-i)w_m\otimes y_0^{(k-p+i)}v_n\\
			&=\sum\limits_{k=p}^{i-1}\,\, {k+r\choose i}\, \,  y(i-k,m-i)w_m\otimes y_0^{(k+r)}v_n.
		\end{array}
	\end{equation}
	
	\noindent Having established \eqref{r>i} and \eqref{r=i-p}, to prove that $\dfrac{W_i}{W_{i-1}}$ is a quotient of $V(\hook(m+n-2i,m-i))$ it suffices to show that \begin{equation}
		\{y(i-l, m-i)w_m \otimes y_{0}^{(r+l)}v_n: 0\leq l\leq i-1\}\subset W_{i-1}, \quad  \text{ for } 1\leq r\leq n-i+1 \label{1.y.r.s}
	\end{equation}
	
	\noindent Since  $w_m\otimes y_0^{(j)}v_n\in W_{i-1}$ for $0\leq j\leq i-1,$ $0\leq i\leq n$,
	$$y(r+i-j,m-i)(w_m\otimes y_0^{(j)}v_n) \in W_{i-1}, \qquad \forall \, \, 0\leq j\leq i-1. $$ Further, for $r\geq i>j\geq 0$ we have,  
	\begin{equation} \label{r>i.j}\begin{array}{ll}
			(-1)^{m-i}x_1^{(m-i)}y_0^{(m+r-j)}(w_m\otimes y_0^{(j)}v_n)&
			= y(i+r-j,m-i)(w_m\otimes y_0^{(j)}v_n),\\ &=\sum\limits_{l=r-j}^{m+r-j-1} {{j+l} \choose {j}} y(r-j+i-l, m-i)w_m \otimes y_{0}^{(j+l)}v_n\\
			&= \sum\limits_{s=0}^{i-1} {{r+s} \choose {j}} y(i-s,m-i)w_m\otimes y_0^{(r+s)}v_n .
	\end{array}\end{equation}
	Setting, $X:= \begin{bmatrix}
		y(i,m-i)w_{m}\otimes y_{0}^{(r)}v_{n}\\
		y(i-1,m_i)w_{m}\otimes y_{0}^{(r+1)}v_{n}\\\\
		\vdots \\
		y(1,m-i)w_{m}\otimes y_{0}^{(r+i-1)}v_{n}\\
	\end{bmatrix}$ and 
	$Y:= \begin{bmatrix}
		(-1)^{m-i}x_1^{(m-i)}y_0^{(m+r)}(w_m\otimes y_0v_n)\\
		(-1)^{m-i}x_1^{(m-i)}y_0^{(m+r-1)}(w_m\otimes y_0^{(1)}v_n)\\\\
		\vdots \\
		(-1)^{m-i}x_1^{(m-i)}y_0^{(m+r-i+1)}(w_m\otimes y_0^{(i-1)}v_n)\\
	\end{bmatrix}$\\
	
	\noindent we can write the set of equations \eqref{r>i.j}, for $0\leq j\leq i-1$, as $Y=A(r,i)X$. 
	As the entries of the colummn matrix $Y$ are contained in $W_{i-1}$ and by \lemref{A(r,i)}, the coeffient matrix $A(r,i)$ is invertible, it follows that 
	the entries of $X$ lie in $W_{i-1}$. Hence, \eqref{1.y.r.s} holds implying that 
	$y(r,m-i)(w_m\otimes y_0^{(i)}v_n)\in U_{i-1}$ for $ i\leq r\leq n-i+1$.\\

	
	\noindent On the other hand, for $r=i-p<i$, $j=r+s \leq i-1$, $1\leq s\leq p-1$, we have 
	\begin{equation} \label{r<j<i}\begin{array}{ll}
			(-1)^{m-i}x_1^{(m-i)}y_0^{(m+r-j)}(w_m\otimes y_0^{(j)}v_n)&= y(i+r-j,m-i)(w_m\otimes y_0^{(j)}v_n),\\ 
			&=\sum\limits_{l=0}^{m+r-j-1} {{r+s+l} \choose {r+s}} y(i-s-l, m-i)w_m \otimes y_{0}^{(r+s+l)}v_n\\
			&=\sum\limits_{k=j-r}^{i-1} {{r+k} \choose {j}} y(i-k, m-i)w_m \otimes y_{0}^{(r+k)}v_n.
	\end{array}\end{equation}
	Now setting, $$\bar{Y}:= \begin{bmatrix}
		(-1)^{m-i}x_1^{(m-i)}y_0^{(m+r)}(w_m\otimes v_n)\\
		\vdots
		\\ (-1)^{m-i}x_1^{(m-i)}y_0^{(m)}(w_m\otimes y_0^{(r)}v_n)\\
		(-1)^{m-i}x_1^{(m-i)}y_0^{(m+r-(r+1)}(w_m\otimes y_0^{(r+1)}v_n)\\
		\vdots \\
		(-1)^{m-i}x_1^{(m-i)}y_0^{(m+r-(i-1)}(w_m\otimes y_0^{(i-1)}v_n)\\
	\end{bmatrix} \qquad \text{and} \qquad \bar{X}:= \begin{bmatrix}
		y(i,m-i)w_{m}\otimes y_{0}^{(r)}v_{n}\\
		y(i-1,m-i)w_{m}\otimes y_{0}^{(r+1)}v_{n}\\\\
		\vdots \\
		y(1,m-i)w_{m}\otimes y_{0}^{(r+i-1)}v_{n}\\
	\end{bmatrix}$$ 
	we can write the equations \eqref{r>i.j} for $0\leq j\leq r$ and \eqref{r<j<i} for $r< j\leq i-1$, as $\bar{Y}=B(r,i)\bar{X}$, where 
	and $B(r,i)$ is the $i\times i$ matrix defined in \lemref{A(r,i)}. As in the case above, using the fact that $B(r,i)$ is invertible (\lemref{A(r,i)}), we see that the entries of 
	$X$ lie in the space spanned by the vectors $\{x_1^{(m-i)}y_0^{(m+r-j)}(w_m\otimes y_0^{(j)}v_n) : 0\leq j\leq i-1\}\subset W_{i-1}.$ 
	
	\noindent If $\overline{ w_m\otimes y_0^{(i)}v_n}$ is the image of $ w_m\otimes y_0^{(i)}v_n$ in $\dfrac{W_i}{W_{i-1}}$ and 
	$\Phi_i:V(\hook(m+n-2i,m-i)\rightarrow \dfrac{W_i}{W_{i-1}}$ is a map such that $\Phi_i(v_{\hook(m+n-2i,m-i)}) = \overline{w_m\otimes y_0^{(i)}v_n},$ then it follows from above that $\Phi_i$ is a surjective $\mathfrak{sl}_2[t]$-module homomorphism.
	Therefore,
	$$\dim W(m\omega)\otimes V(n) = \sum\limits_{i=0}^n \dim \dfrac{W_i}{W_{i-1}} \leq \sum\limits_{i=0}^n \dim V(\hook(m+n-2i,m-i)).$$ 
	Setting $W_{-1}=0$ and using \lemref{gen.W_loc.x.V(p)} and \lemref{sum dim Ui/Ui-1}(i), we thus see that $$\sum\limits_{i=0}^n \dim \dfrac{W_i}{W_{i-1}}= 2^m(n+1).$$ However, this is possible only if the surjection $\Phi_i$ is an isomorphism for $0\leq i\leq n$. This completes the proof of the theorem in this case.\\
	
	When $m\leq n$, considering the generator $y_{0}^{m} w_{m}\otimes v_{n}$ of {$W_{loc}(m) \otimes V(n) $} and setting $U_{i} = U(\mathfrak{sl}_{2}[t]) (y_{0}^{(m-i)}w_m \otimes v_n)$ for $0\leq i\leq m$, we see that there exists a decreasing chain 
	of submodules $$ W_{loc}(m \omega)\otimes V(n)= U_{0} \supset U_{1} \supset \cdots \supset U_{m} \supset 0.$$ Using similar arguments as above and  \lemref{sum dim Ui/Ui-1}(ii) we see that even in this case, $\dfrac{U_{i}}{U_{i+1}}$ is isomorphic to $V(\hook(n-m+2i, i))$ for $ 1 \leq i \leq m$ and $\dfrac{U_{0}}{U_{1}}$ is isomorphic to  $V(n-m)$. This completes the proof of the theorem. \endproof
	
	\vspace{.25cm}
	
	\noindent The following result is an immediate consequence of \thmref{filtration in irreducible}
	
	\begin{corollary} \label{Character W_m.V_n} Given two positive integers $m,n$, we have,
		$$\begin{array}{l}ch_{gr}\, W_{loc}(m) \otimes V(n) \\ 
			\, \qquad = \left\{ \begin{array}{ll}\sum\limits_{i=0}^n ch_{gr} \, 
				V(i+1, 1^{m-n+i-1}) & \text{for } m> n, \\
				\sum\limits_{i=1}^{m} ch_{gr} \, V(n-m+i+1, 1^{i-1}) + ch_{gr}\, V(n-m, 0) & \text{for } m\leq n.
			\end{array} \right.\end{array} $$
	\end{corollary}
	
	\begin{remark} Combinatorial formulas for determining the multiplicity of level $\ell$ Demazure modules in CV-modules associated with hook-type partitions have been obtained in \cite[Theorem 4]{MR4229660} and \cite[Theorem 4.3]{MR3210603}. Using the latter in combination with \corref{Character W_m.V_n}, one can deduce a combinatorial formula for the graded multiplicities of level $\ell$ Demazure modules in $W_{loc}(m\omega)\otimes ev_{0}V(n\omega)$ for all positive integers $m,n$. \end{remark}
	
	\subsection{} In view of \thmref{present},  we shall refer to $W_{loc}(m)\otimes V(n)$ as $W(m,n)$ in the rest of this section. As $W(m,n)$ is a graded $\mathfrak{sl}_2[t]$-module, 
	$$ch_{gr} W(m,n) = \sum\limits_{k\in \mathbb Z_+} ch_{\mathfrak{sl}_2}\ W(m,n)[k]\ q^k,$$
	where $W(m,n)[k]$ is a finite-dimensional $\mathfrak{sl}_2$-module.\\ 
	Define $[W(m,n):V(r)]_q$ as the polynomial in indeterminate $q$ given by
	$$[W(m,n):V(r)]_q = \sum\limits_{k\geq 0}\, [W(m,n)[k]: V(r,k)] q^k, $$
	where $[W(m,n)[k]: V(r,k)]$, denotes  the multiplicity of the irreducible $\mathfrak{sl}_2[t]$-module $V(r,k)$ in $W(m,n)[k]$ for $r\in \mathbb Z_+$. Using \thmref{filtration in irreducible}, we now obtain the polynomials $[W(m,n):V(r)]_q$.
	
	\begin{theorem}\label{graded mul of irr in W(m,n)}  Let $m,n\in \mathbb N$,  and $i\in \mathbb Z_+$ be such that $m+n-2i\geq 0$.\\
		If $m \leq n$,	 
		$$[W(m,n): V(m+n-2i)]_q =  \begin{bmatrix} m \\ i \end{bmatrix}_{q}, \quad 0 \leq i \leq m ,$$
		If $m>n$,
		$$\begin{array}{l}[W(m,n): V(m+n-2i)]_q =\left\{\begin{array}{ll}  \begin{bmatrix} m \\ i \end{bmatrix}_{q}, \quad 0 \leq i \leq n,\\
				\begin{bmatrix} m \\ i \end{bmatrix}_{q} -\begin{bmatrix} m \\ i-n-1 \end{bmatrix}_{q}, \quad  n+1 \leq i \leq \lfloor{\frac{n+m}{2}} \rfloor\end{array}\right.  \end{array}$$
	\end{theorem}
	The proof of the theorem occupies the rest of this section.\\
	
	\vspace{.25cm}
	\noindent At the onset we record {three} $q-$identities that we shall use. 
	\begin{eqnarray}
		\label{q-combn identity.3} 
		&q^r\begin{bmatrix} n-1 \\ r \end{bmatrix}_{q}+\begin{bmatrix} n-1 \\ r-1 \end{bmatrix}_{q}= \begin{bmatrix} n \\ r \end{bmatrix}_{q}, &\qquad \text{ for } n\in \mathbb N, \, r\in \mathbb Z_+.\\
		\label{q-combn identity}
		& \qquad \begin{bmatrix} n \\ r \end{bmatrix}_{q}+q^{n-r+1}\begin{bmatrix} n \\ r-1 \end{bmatrix}_{q}= \begin{bmatrix} n+1 \\ r \end{bmatrix}_{q}, &\qquad \text{ for } n\in \mathbb N, \, r\in \mathbb Z_+.   \\
		&\sum\limits_{i=0}^{n} q^{i}  \begin{bmatrix} m+i \\ i \end{bmatrix}_{q} =  \begin{bmatrix} n+m+1 \\ n \end{bmatrix}_{q} &\qquad \text{ for } m\in \mathbb N, \, n\in \mathbb Z_+. \label{q-combn identity.2}
	\end{eqnarray}
	
	\noindent We need the following results to prove \thmref{graded mul of irr in W(m,n)}. 
	\begin{lemma} \label{Kus thm} Given $k \in \mathbb{N} \text{ and } r \in Z_{+}$, we have, 
		$$\begin{array}{lll}
			ch_{gr} V(k+r,1^k) &= & ch_{gr} V(k+r+1, 1^{k-1}) + q^k \Big( ch_{\mathfrak{sl}_2} V(r,0) + 
			\sum\limits_{i=2}^k ch_{gr} V(i+r,1^{i-2})\Big)\\ &\\
			ch_{gr} V(k,1^{k+r}) &= & ch_{gr}\ V(k+1, 1^{k+r-1}) +q^{k+r}  \Big( ch_{gr} W_{loc}(r)+ \sum\limits_{i=2}^k ch_{gr}\ V(i,1^{i+r-2})\Big)\end{array} $$
	\end{lemma}
	\proof This is a direct consequence of \cite[Theorem 18]{MR4229660}.\endproof
	
	\begin{lemma} \label{character of hooks in terms of irreducibles} 
		The graded characters of $V(\xi)$ corresponding to hook type partitions 
		$\xi$ are as follows:
		$$\begin{array}{lll} ch_{gr} V(k+r,1^k) &= & 
			\sum\limits_{p=0}^{k} q^p \begin{bmatrix} k \\ p \end{bmatrix}_{q} ch_{\mathfrak{sl}_2} V(r+2(k-p))\\ & & \\
			ch_{gr} V(k,1^{k+1}) & = & \sum\limits_{p=0}^{k} q^{k-p} \begin{bmatrix} k+1 \\ p+1 \end{bmatrix}_{q} ch_{\mathfrak{sl}_2} V(1+2p)\\ & &\\
			ch_{gr} V(k,1^{k+r}) & = & \sum\limits_{p=0}^{k} q^p \begin{bmatrix} k+r \\ p \end{bmatrix}_{q} ch_{\mathfrak{sl}_2} V(r+2(k-p))\\
			&  +& \sum\limits_{p=k+1}^{\lfloor\frac{2k+r}{2}\rfloor} q^p (\begin{bmatrix} k+r \\ p \end{bmatrix}_{q} - q^{2k+r+1-2p}\begin{bmatrix} k+r \\ p-k-1 \end{bmatrix}) ch_{\mathfrak{sl}_2}  V(r-2(p-k)) \qquad \text{ for } r\geq 2.
		\end{array} $$ \end{lemma} 
	\proof Applying induction on $k$, we prove the character formula of $V(k+r,1^k)$,  $r\geq 0$.\\
	When $k=1$, by \lemref{CV iso}(b) we have the following short exact 
	sequence : 
	$$0 \rightarrow \tau_{1}(V(n-1)) \rightarrow V(n,1) \xrightarrow{\phi} V(n+1,0)\rightarrow 0$$
	Therefore, $$ch_{gr} V(n,1) = ch_{\mathfrak{sl}_2} V(n+1) + q\ ch_{\mathfrak{sl}_2} V(n-1) = \sum\limits_{p=0}^{1}q^p \begin{bmatrix} 1 \\ p \end{bmatrix}_{q} ch_{\mathfrak{sl}_2}\ V(n+1-2p).$$
	Assume that the given character formula holds for all $k<l$ and consider the module $V(l+r,l)$. Using \lemref{Kus thm} and applying the inductive 
	hypothesis, we get \\
	\begin{align*}
		& ch_{gr}\ V(l+r,1^l) = ch_{gr} V(l+r+1,1^{l-1}) + q^l ch_{\mathfrak{sl}_2} V(r) + 
		q^l \big(\sum\limits_{i=2}^{l} ch_{gr} V(i-2+r+2,1^{i-2})\big) \\
		&= \sum\limits_{p=0}^{l-1} q^p \begin{bmatrix} l-1 \\ p \end{bmatrix}_{q} ch_{\mathfrak{sl}_2} V(r+2+2(l-p))  + q^l \Big(ch_{\mathfrak{sl}_2} V(r) +\sum\limits_{i=0}^{l-2} 
		\big(\sum\limits_{s=0}^{i} q^s \begin{bmatrix} i \\ s \end{bmatrix}_{q} 
		ch_{\mathfrak{sl}_2}\ V(r+2+2i-2s)  \big) \Big) \\
		& = \sum\limits_{p=0}^{l-1} q^p \begin{bmatrix} l-1 \\ p \end{bmatrix}_{q} ch_{\mathfrak{sl}_2} V(r+2(l-p+1)) + q^l \Big(  ch_{\mathfrak{sl}_2} V(r)+   \sum\limits_{j=0}^{l-2} \big(\sum\limits_{s=0}^{l-2-j} q^s \begin{bmatrix} j+s \\ s \end{bmatrix}_{q} ch_{\mathfrak{sl}_2}\ V(r+2+2j)  \big) \Big) \\ 
		& =  \sum\limits_{j=1}^{l} q^{l-j} \begin{bmatrix} l-1 \\ l-j \end{bmatrix}_{q} ch_{\mathfrak{sl}_2} V(r+2(j+1))  + \sum\limits_{j=0}^{l-2} q^l  \begin{bmatrix} l-1 \\ j+1 \end{bmatrix}_{q} 
		ch_{\mathfrak{sl}_2}\ V(r+2(j+1)) + q^lch_{\mathfrak{sl}_2} V(r), \quad \text{ by } \eqref{q-combn identity.2} \\ 
		& = \sum_{p=0}^{l} q^p \begin{bmatrix} l \\ p \end{bmatrix}_{q} ch_{\mathfrak{sl}_2} 
		V(r+2(l-p)) \qquad \quad \text{ by } \eqref{q-combn identity}. \end{align*}
	This establishes the character formula for the module $V(k+r,1^k)$.\\
	
	\noindent When $r=1$, using \lemref{Kus thm} and the character formula for $V(k+1,1^k)$, we have
	\begin{align*}
		&ch_{gr}(V(k,1^{k+1} )) = ch_{gr} V(k+1,1^{k}) +  
		q^{k+1} \big(\sum\limits_{i=1}^{k} ch_{gr} V(i,1^{i-1})\big) \\
		& = \sum\limits_{p=0}^{k} q^p \begin{bmatrix} k \\ p \end{bmatrix}_{q} ch_{\mathfrak{sl}_2} V(1+2(k-p))  + q^{k+1} \Big(\sum\limits_{j=0}^{k-1} 
		\big(\sum\limits_{s=0}^{j} q^s 
		\begin{bmatrix} j \\ s \end{bmatrix}_{q} ch_{\mathfrak{sl}_2} V(1+2(j-s)) \big) \Big) \\&  (\text{ Rearranging the terms and using \eqref{q-combn identity.2} we get })\\
		& = \sum\limits_{p=0}^{k} q^p \begin{bmatrix} k \\ p \end{bmatrix}_{q} ch_{\mathfrak{sl}_2} V(1+2(k-p))  + q^{k+1} \Big(\sum\limits_{i=0}^{k-1}  
		\begin{bmatrix} k \\ i+1 \end{bmatrix}_{q} ch_{\mathfrak{sl}_2} V(1+2i) \Big)
		\\& = \sum\limits_{j=0}^k q^{k-j} \Big(\begin{bmatrix} k \\ j \end{bmatrix}_{q} + q^{j+1} \begin{bmatrix} k \\ j+1 \end{bmatrix}_{q}\Big) ch_{\mathfrak{sl}_2} V(1+2j) = \sum\limits_{j=0}^k q^{k-j} \begin{bmatrix} k+1 \\ j+1 \end{bmatrix}_{q} V(1+2j),  \qquad (\text{using \eqref{q-combn identity.3}}).
	\end{align*} 
	To to establish the character formula of $V(k,1^{k+r})$ for $r\geq 2$, 
	we apply induction on $r$.\\
	When $r=2$, using \lemref{Kus thm} and the character formula for $V(k,1^k)$ we have, 
	\begin{align*}
		&ch_{gr}(V(k,1^{k+2} )) = ch_{gr} V(k+1,1^{k+1}) + q^{k+2} \big(\sum\limits_{i=2}^{k} ch_{gr} V(i,1^{i-1})\big) \\
		&=\sum\limits_{p=0}^{k+1}  \begin{bmatrix} k+1 \\ p \end{bmatrix}_{q} ch_{\mathfrak{sl}_2} V(2(k+1-p)) + q^{k+2}\Big( \sum\limits_{i=1}^k \sum\limits_{s=0}^i \begin{bmatrix} i \\ p \end{bmatrix}_{q} ch_{\mathfrak{sl}_2} V(2(i-p)) \Big)
	\end{align*}
	Rearranging the terms in the summation and using \eqref{q-combn identity.2}
	we have,
	\begin{align*}
		&ch_{gr}(V(k,1^{k+2} ))=\sum\limits_{p=0}^k q^p \begin{bmatrix} k+1 \\ p \end{bmatrix}_{q} ch_{\mathfrak{sl}_2} V(2+2(k-p)) \\&+ \sum\limits_{j=1}^{k} q^{k+2}\begin{bmatrix} k+1 \\ j+1 \end{bmatrix}_{q} ch_{\mathfrak{sl}_2} V(2j) + (q^{k+1} +q^{k+3})\begin{bmatrix} k \\ 1 \end{bmatrix}_{q} ch_{\mathfrak{sl}_2} V(0)\\
		&= V(2k+2) + (q^{k+1} +q^{k+3})\begin{bmatrix} k \\ 1 \end{bmatrix}_{q} ch_{\mathfrak{sl}_2} V(0)+\sum\limits_{j=1}^k q^{k-j+1} \big( \begin{bmatrix} k+1 \\ j \end{bmatrix}_{q} + q^{j+1}\begin{bmatrix} k+1 \\ j+1 \end{bmatrix}_{q}\big) ch_{\mathfrak{sl}_2} V(2j) \\
		\\
		& \text{Using \eqref{q-combn identity.3} and substituting $k-j+1=p$ we get, }\\
		&=  \sum\limits_{j=1}^k q^{p} \begin{bmatrix} k+2 \\ k-p+2 \end{bmatrix}_{q}
		ch_{\mathfrak{sl}_2} V(2+2(k-p)) + q^{k+1}\big( \begin{bmatrix} k+2 \\ 1 \end{bmatrix}_{q} -q \big) ch_{\mathfrak{sl}_2} V(0).
	\end{align*} 
	This proves that the formula for $V(k,1^{k+r})$ holds when $r=2$. Similarly, we can prove the formula for $r=3$ separately.
	Now, assume that the character formula for $V(k,1^{k+r})$ holds for all $4 \leq r<l$. Then
		\begin{align*}
			&ch_{gr} V(k,1^{k+l}) =  ch_{gr} V(k+1, 1^{k+l-1}) +q^{k+l}  \Big( \sum\limits_{i=1}^k ch_{gr} V(i,1^{i+l-2})\Big) \\
			&= \sum\limits_{p=0}^{k+1} q^p 
			\begin{bmatrix} k+l-1 \\ p \end{bmatrix}_{q} ch_{\mathfrak{sl}_2} V(2k+l-2p) + 
			q^{k+l}
			\Big(\sum\limits_{i=1}^k \big( \sum\limits_{p=0}^{i} q^p \begin{bmatrix} i+l-2 \\ p \end{bmatrix}_{q} ch_{\mathfrak{sl}_2} V(l+2(i-p-1))\Big)\\ 
			&  + \sum\limits_{p=k+2}^{\lfloor\frac{2k+l}{2}\rfloor} q^p (\begin{bmatrix} k+l-1 \\ p \end{bmatrix}_{q} - q^{2k+l+1-2p}\begin{bmatrix} k+l-1 \\ p-k \end{bmatrix}) 
			ch_{\mathfrak{sl}_2}  V(2k+l-2p)  \\ 
			& + q^{k+l}\Big(\sum\limits_{p=i+1}^{\lfloor\frac{2i+l-2}{2}\rfloor} q^p (\begin{bmatrix} i+l-2 \\ p \end{bmatrix}_{q} - q^{2i+l-1-2p}\begin{bmatrix} i+l-2 \\ p-i-1 \end{bmatrix}) ch_{\mathfrak{sl}_2}  V(2i+l-2-2p)   \big)\Big)  \\
			&=\sum\limits_{p=0}^{k} q^{p} \begin{bmatrix} k+l \\ p \end{bmatrix}_{q} ch_{gr}V(2k+l-2p,0) \\&+ \sum_{p= k+1}^{\lfloor {\frac{2k+l}{2}} \rfloor} q^{p} (\begin{bmatrix} k+l \\ p \end{bmatrix}_{q} - q^{2k+l+1-2p} \begin{bmatrix} k+l \\ p-k-1 \end{bmatrix}_{q}) ch_{gr}V(2k+l-2p,0),
		\end{align*}
		where the last line follows after rearranging the terms and using \eqref{q-combn identity.3} and \eqref{q-combn identity}. \endproof
		
		We now complete the proof of \thmref{graded mul of irr in W(m,n)}.
		\proof We first consider the case when $m\leq n$. Suppose that $n=m+r$.\\ 
		Then from \corref{Character W_m.V_n} and \lemref{character of hooks in terms of irreducibles} it follows that :
		\begin{align*}
			ch_{gr}\ W(m,m+r) & = \sum\limits_{i=1}^{m} ch_{gr} V(r+2+i-1,1^{i-1})+ ch_{\mathfrak{sl}_2}\ V(r)\\
			&= \sum\limits_{k=0}^{m-1}\Big( \sum\limits_{p=0}^{k} q^{p} \begin{bmatrix} k \\ p \end{bmatrix}_{q} ch_{\mathfrak{sl}_2}\ V(r+2+2(k-p))\Big)+ ch_{\mathfrak{sl}_2}V(r) \\
			\text{Rearranging the terms we thus get},\\
			& = \sum\limits_{i=1}^{m} \Big(\sum\limits_{p=0}^{m-i} q^p \begin{bmatrix} p+i-1 \\ p \end{bmatrix}_{q}\Big) ch_{\mathfrak{sl}_2}\ V(r+2+2i) + ch_{\mathfrak{sl}_2} V(r)\\
			&=\sum_{i=0}^{m} \begin{bmatrix} m \\ i \end{bmatrix}_{q} Ch_{gr}V(r+2i), \quad (\text{ by } \eqref{q-combn identity.2} ) .
		\end{align*}
		
		
		\noindent We now consider the case when $m> n$, i.e., $m=n+r$.  It follows from \corref{Character W_m.V_n} that,  
		$$	ch_{gr} W(n+r,n) = \sum\limits_{i=0}^{n} ch_{gr} V(i+1,1^{r+i-1} )= \sum\limits_{k=1}^{n+1} ch_{gr} V(k,1^{r-2+k}).$$
		
		\noindent Hence, using  \lemref{character of hooks in terms of irreducibles}, for $r>3$,  we have,
		$$\begin{array}{l}
			ch_{gr} W(n+r,n) 
			= \sum\limits_{k=1}^{n+1}\big(\sum\limits_{p=0}^k q^p \begin{bmatrix} k+r-2 \\ p \end{bmatrix}_{q} ch_{\mathfrak{sl}_2} V(r-2+2(k-p)) \big)\\
			\quad + \sum\limits_{k=1}^{n+1}\Big( \sum\limits_{p=k+1}^{\lfloor\frac{2k+r-2}{2}\rfloor} q^p \big( \begin{bmatrix} k+r-2 \\ p \end{bmatrix}_{q} - q^{r-1-2(p-k)}\begin{bmatrix} k+r-2 \\ p-k-1 \end{bmatrix} \big) ch_{\mathfrak{sl}_2} V(r-2-2(p-k))\Big)\\ 
			= \sum\limits_{i=-1}^n \big(\sum\limits_{p=0}^{n-i} q^p \begin{bmatrix} p+i+r-1 \\ p \end{bmatrix}_{q} \big) ch_{\mathfrak{sl}_2} V(r+2i)\\
			+ \sum\limits_{i=1}^{\lfloor\frac{r}{2} \rfloor-1} \sum\limits_{p=i+1}^{i+n+1} 
			q^p\big(\begin{bmatrix} p-i+r-2 \\ p \end{bmatrix}_{q} -q^{r-1-2i}\begin{bmatrix} p-i+r-2 \\ i-1 \end{bmatrix}_{q} \big) ch_{\mathfrak{sl}_2} V(r-2(i+1))\\
			= \sum\limits_{i=-1}^n\begin{bmatrix} n+r \\ n-i \end{bmatrix}_{q}  ch_{\mathfrak{sl}_2} V(r+2i) \\ + \sum\limits_{i=1}^{\lfloor\frac{r}{2} \rfloor-1} \Big( \begin{bmatrix} n+r \\ n+i+1 \end{bmatrix}_{q} - \begin{bmatrix} r-1 \\ i \end{bmatrix}_{q}  - \big( \sum\limits_{j=0}^n q^{r-i+j} \begin{bmatrix} r-1+j \\ i-1 \end{bmatrix}_{q}\big) \Big) ch_{\mathfrak{sl}_2} V(r-2(i+1)) \quad (\text{ by } \eqref{q-combn identity.2})\\
			= \sum\limits_{i=-1}^n\begin{bmatrix} n+r \\ n-i \end{bmatrix}_{q}  ch_{\mathfrak{sl}_2} V(r+2i)  + \sum\limits_{i=1}^{\lfloor\frac{r}{2} \rfloor-1} \Big( \begin{bmatrix} n+r \\ n+i+1 \end{bmatrix}_{q} - \begin{bmatrix} r+n \\ i-1 \end{bmatrix}_{q}  \Big) ch_{\mathfrak{sl}_2} V(r-2(i+1)) \quad (\text{ by } \eqref{q-combn identity})
		\end{array}$$
		The proof of the theorem in the case when $r=3$, is similar except for the fact that here we need to use the formula for $ch_{gr} V(k,1^{k+1})$ as given in \lemref{character of hooks in terms of irreducibles}. 
		
		For $m=n+1$, using \corref{Character W_m.V_n}, \lemref{character of hooks in terms of irreducibles} and identities \eqref{q-combn identity.3}-\eqref{q-combn identity.2},  we have:
		\begin{align*}
			&ch_{gr} W(n+1,n) = \sum\limits_{i=0}^n ch_{gr} V(i+1,1^i) = \sum\limits_{i=0}^n \sum\limits_{p=0}^i q^p \begin{bmatrix} i \\ p \end{bmatrix}_{q} ch_{\mathfrak{sl}_2} V(1+2(i-p))\\
			&=\sum\limits_{j=0}^n \sum\limits_{p=0}^{n-j} q^p \begin{bmatrix} p+j \\ p \end{bmatrix}_{q} ch_{\mathfrak{sl}_2} V(1+2j) = \sum\limits_{j=0}^{n} \begin{bmatrix} n+1 \\ n-j \end{bmatrix}_{q} ch_{\mathfrak{sl}_2} V(1+2j).
		\end{align*}
		Likewise, for $m=n+2$, we have, 
		\begin{align*}
			&ch_{gr} W(n+2,n) = \sum\limits_{i=0}^n ch_{gr} V(i+1,1^{i+1}) = \sum\limits_{i=0}^n \sum\limits_{p=0}^{i+1} q^p \begin{bmatrix} i+1 \\ p \end{bmatrix}_{q} ch_{\mathfrak{sl}_2} V(2+2(i-p))\\
			&=\sum\limits_{j=0}^n \sum\limits_{p=0}^{n-j} q^p \begin{bmatrix} p+j+1 \\ p \end{bmatrix}_{q} ch_{\mathfrak{sl}_2} V(2+2j) + \sum\limits_{p=1}^{n+1} q^p \begin{bmatrix} p\\ p \end{bmatrix}_{q} ch_{\mathfrak{sl}_2} V(0)\\
			&= \sum\limits_{j=0}^{n} \begin{bmatrix} n+2 \\ n-j \end{bmatrix}_{q} ch_{\mathfrak{sl}_2} V(2+2j) + \big(\begin{bmatrix} n+2\\ 1 \end{bmatrix}_{q}-1\big) ch_{\mathfrak{sl}_2} V(0) .
		\end{align*}
		This completes the proof of \thmref{graded mul of irr in W(m,n)}.\endproof

		
		\section{The module $W_{loc}(n\omega)\otimes W_{loc}(m\omega)$ and Demazure flags}
		
		In this section we recall the definition of truncated local Weyl modules and show that the tensor product of two local Weyl modules has a filtration by truncated local Weyl modules. Then using \propref{Theorem 3.3}, we establish these tensor product modules have an excellent filtration and give a closed formula for the graded multiplicity of level two Demazure modules in the tensor product of level one Demazure modules. In addition we show that the tensor product of a level one Demazure module with a level two Demazure module admits a filtration by level 3 Demazure modules.
		
		\subsection{Truncated local Weyl Modules}\label{truncated} We begin by recalling the definition of truncated local Weyl modules.\\ Given a pair of integers $(m,N)$, the truncated local Weyl module $W_{loc}([m],N)$ is a quotient of the local Weyl module $W_{loc}(m)$ generated by an element $w_{m,N}$ satisfying the following relations:
		\begin{equation}\label{Trunc.eq}
			\begin{array}{c}(x\otimes t^r)w_{m,N}=0, \qquad (h\otimes t^r)w_{m,N}=m\delta_{r,0}w_{m,N}, \qquad (y\otimes 1)^{m+1}w_{m,N} =0,\\
				y\otimes t^s.w_{m,N}=0, \quad \forall\, s\geq N. \end{array}
		\end{equation}
		The truncated local Weyl modules were studied in \cite[Theorem 4.3]{MR3407180} and it was shown that if $m=Nk+r$ with $0\leq r<N$,  then 
		the $\mathfrak{sl}_2[t]$-module $W_{loc}([m],N)$ is isomorphic to the fusion product module $V(k+1)^{\ast r}\ast V(k)^{\ast N-r}$. Hence by \cite[Theorem 5]{MR3296163}(ii), it is known that the module $W_{loc}([m],N)$ is isomorphic to the CV-module $V(\xi(m,N))$ where $\xi(m,N)=((k+1)^r, k^{N-r}).$
		
		\vspace{.25cm}
		In particular, given non-negative integers $a$ and $b$, since $2a + b = 1 \cdot (a + b) + a$, the module $W_{\text{loc}}([2a + b], a + b)$ is isomorphic to $V(2^a, 1^b)$ and hence has a Demazure flag of level $2$, as per \propref{Theorem 3.3}. The following result gives the  character formula for such truncated local Weyl modules in terms of level 1 and level 2 Demazure modules.
		
		\begin{lemma} \label{char V(a,b)} Give two non-negative integers $a,b$, we have,\\
			i. $ch_{gr} V(2^a,1^b)= \sum\limits_{k=0}^{a} (-1)^{k} \begin{bmatrix} a \\ k \end{bmatrix}_q q^{k(a+b)-\frac{k(k-1)}{2}}\ ch_{gr}\ W_{loc}(b+2a-2k)$\\
			
			\noindent ii. The truncated local Weyl module $W_{loc}([2a+b],a+b)$ has a level two Demazure flag and 
			\begin{equation}
				ch_{gr}V(2^a,1^b) = \sum\limits_{k=0}^{\lfloor{\frac{b}{2}}\rfloor}   q^{k(a+\lceil{\frac{b}{2}}\rceil)}  \begin{bmatrix} {\lfloor{\frac{b}{2}}\rfloor} \\ k \end{bmatrix}_q \, ch_{gr} D(2,2a+b-2k)
			\end{equation}
		\end{lemma}
		\proof i. We prove part (i) of  the lemma by applying induction on $a$.\\
		For $a=1$, from the short exact sequence of $W_{loc}((b+2)\omega)$ given in \lemref{CV iso}(b), we deduce that 
		$$ch_{gr}\ V(2,1^{b}) = ch_{gr} W_{loc}(b+2) - q^{(b+1)} ch_{gr} W_{loc}(b) $$
		Thus, the result holds for $a=1$. \\
		Assume that the result holds when $a<l$ and consider the case when $a=l$.\\
		Once again, from \lemref{CV iso}(b) we deduce that 
			\begin{align*}	
				ch_{gr} V(2^l , 1^{b}) & = ch_{gr}\ V(2^{l-1},1^{b+2}) - q^{(l+b)} ch_{gr}\ V(2^{l-1}, 1^b)\\
				&\text{Applying induction hypothesis we have,}\\
				& = \sum\limits_{k=0}^{l-1} (-1)^{k} \begin{bmatrix} l-1 \\ k \end{bmatrix}_q q^{k(l+b+1)- \frac{k(k-1)}{2}} ch_{gr}\ W_{loc}(b+2(l-k))\\
				&- q^{l+b} \sum\limits_{k=0}^{l-1} (-1)^{k} \begin{bmatrix} l-1 \\ k \end{bmatrix}_q q^{k(l+b-1)- \frac{k(k-1)}{2}}ch_{gr} W_{loc}(b+2(l-k-1)) \\
				& = \sum\limits_{k=0}^{l-1} (-1)^{k} \begin{bmatrix} l-1 \\ k \end{bmatrix}_q q^{k(l+b+1)- \frac{k(k-1)}{2}} ch_{gr}\ W_{loc}(b+2(l-k))\\
				& - \sum\limits_{k=0}^{l-1} (-1)^{k} \begin{bmatrix} l-1 \\ k \end{bmatrix}_q q^{(k+1)(l+b)- \frac{k(k+1)}{2}} ch_{gr}\ W_{loc}(b+2(-k-1))
			\end{align*}
			\begin{align*}
				\text{Reindexing  the indices we get},\\
				ch_{gr}V(2^l , 1^b )	& = \sum\limits_{k=0}^{l-1} (-1)^{k} \begin{bmatrix} l-1 \\ k \end{bmatrix}_q q^{k(l+b+1)- \frac{k(k-1)}{2}} ch_{gr} W_{loc}(b+2(l-k))\\
				& + \sum\limits_{k=1}^{l} (-1)^{k} \begin{bmatrix} l-1 \\ k-1 \end{bmatrix}_q q^{k(l+b)- \frac{k(k-1)}{2}} ch_{gr} W_{loc}(b+2(l-k)) \\
				& = \sum\limits_{k=0}^{l} (-1)^{k} \begin{bmatrix} l
					\\ k \end{bmatrix}_q q^{k(l+b)- \frac{k(k-1)}{2}} ch_{gr} W_{loc}(b+2(l-k)). \qquad \qquad 	\end{align*}
			
			\noindent ii. We prove part(ii) of the lemma by applying induction on $b$.\\
			First note that for any partition $\bon$ of $n\in \mathbb N$, the multiplicity of a Demazure module $D(m,r)$ with highest weight $r$ is zero whenever $n-r\notin 2\mathbb Z_+,$ and $$V(2^a) \cong D(2,2a), \qquad V(2^a,1)\cong D(2,2a+1).$$ Hence there is nothing to prove for $b=0$ and 1.\\
			Now suppose $2a+b\in 2\mathbb Z+1$ and $b=2d+1$. \\
			For $d=1$,  using (\ref{CV iso}), we have the following short exact sequence,
			$$0 \rightarrow \tau_{a+1}V(2^{a},1) \rightarrow V(2^a, 1^{3}) \xrightarrow{\phi} V(2^{a+1},1)\rightarrow 0 $$
			Hence,
			$$\begin{array}{ll}ch_{gr} V(2^a, 1^{3} ) &= ch_{gr}D(2, 2(a+1)+1) + q^{a+d+1}ch_{gr}D(2,2a+1) \\&= \sum\limits_{k=0}^{1} q^{k(a+d+1)} \begin{bmatrix} d \\ k \end{bmatrix}_{q}ch_{gr}D(2,2a+3-2k)\end{array}$$
			\\
			Assume that the result holds for all $d<l$ and consider module $V(2^a,1^{2l+1})$.\\
			Once again, using (\ref{CV iso}), we have the following short exact sequence,
			$$0 \rightarrow \tau_{a+2b-1}V(2^{a},1^{2l-1}) \rightarrow V(2^a, 1^{2l+1}) \xrightarrow{\phi} V(2^{a+1},1^{2l-1})\rightarrow 0$$
			\begin{align*}
				ch_{gr}V(2^a,1^{2l+1}) & = ch_{gr}V(2^{a+1}, 1^{2l-1}) + q^{a+2l} ch_{gr}V(2^a, 1^{2l-1})\\
				& = \sum\limits_{k=0}^{l-1}  q^{k(a+l+1)} \begin{bmatrix} l-1 \\ k \end{bmatrix}_{q} ch_{gr} D(2,1+2(a+l-k)) \quad \text{(By induction hypothesis)}\\
				& +  \sum_{r=0}^{l-1} q^{r(a+l)+a+2l} \begin{bmatrix} l-1 \\ r \end{bmatrix}_{q} ch_{gr} D(2,1+2(a+l-k-1))\\ 
				& = \sum_{k=0}^{l}  q^{k(a+\lceil\frac{b}{2}\rceil)} \begin{bmatrix} \lfloor\frac{b}{2}\rfloor \\ k \end{bmatrix}_{q} ch_{gr} D(2,2a+b-2k),
			\end{align*} 
			where we obtain the last equality by using the q-combinatorial 
			identity \eqref{q-combn identity}. This proves part (ii) when $2a + b$ is an odd integer. The proof is similar in the case when $2a + b$ is an even integer.
			\endproof
			
			
			\vspace{.15cm}
			
			\subsection{Character formula of $W_{loc}(m)\otimes W_{loc}(n)$.}The following result on the character formula for the tensor product of two local Weyl modules of $\mathfrak{sl}_2[t]$ was given in \cite{9}. Since this result plays an important role in the proof of the main result of this section, \thmref{tensor filtration}, for the sake of completeness, we include its independent proof here. 
			
			\begin{lemma} \label{tensor of two local} Given $m,n\in \mathbb Z_+$,
				\begin{equation} \label{tensor p.rule}
					ch_{gr} W_{loc}(n)\otimes W_{loc}(m) =
					\sum_{i=0}^{\min\{n,m\}} \begin{bmatrix} n \\ i \end{bmatrix}_q \begin{bmatrix} m \\ i \end{bmatrix}_q (1-q)\dots (1-q^i ) ch_{gr} W_{loc}(n+m-2i).
				\end{equation}
			\end{lemma}
			\proof The result is a consequence of \cite[Proposition 3.2, Proposition 5.13]{MR3391927}.
			
			\begin{remark} In Appendix we give an alternate proof of the lemma using Pieri rule.
			\end{remark}
			\begin{lemma}\label{Blanton.result}
				For $m,n \in \mathbb{Z}_{+}$ with  $n\geq m$, 
				\begin{equation}\label{local tensor char}
					ch_{gr}(W_{loc}(m) \otimes W_{loc}(n)) = \sum_{k=0}^{m} \begin{bmatrix} m \\ k \end{bmatrix}_{q} ch_{gr} V(2^{m-k},1^{n-m}). \end{equation}
			\end{lemma}
			\proof Given \eqref{tensor p.rule}, to prove \eqref{local tensor char} it suffices to show that
			$$\begin{array}{l} \sum\limits_{i=0}^{m} \begin{bmatrix} n \\ i \end{bmatrix}_q \begin{bmatrix} m \\ i \end{bmatrix}_q (1-q)\dots (1-q^i )\, ch_{gr}W_{loc}(n+m-2i)=\sum\limits_{k=0}^{m} \begin{bmatrix} m \\ k \end{bmatrix}_{q} ch_{gr} V(2^{m-k},1^{n-k}), \end{array}$$ or equivalently: 
			
			\begin{equation} \begin{array}{l} \sum\limits_{i=0}^{m} \begin{bmatrix} n \\ i \end{bmatrix}_q \begin{bmatrix} m \\ i \end{bmatrix}_q (1-q)\dots (1-q^i )\, ch_{gr}W_{loc}(n+m-2i)\\ = \sum\limits_{k=0}^{m} \sum\limits_{i=0}^{k} (-1)^{i} \begin{bmatrix} m \\ k \end{bmatrix}_q \begin{bmatrix} k \\ i \end{bmatrix}_q q^{i(k+n-m)- \frac{i(i-1)}{2}}\, ch_{gr}W_{loc}(n-m+2k-2i)  \end{array} \label{tensor q eq}\end{equation} Using \lemref{char V(a,b)}. Rearranging the terms on the  right hand side of \eqref{tensor q eq} we have:
			
			$$\begin{array}{l}
				\sum\limits_{k=0}^{m} \sum\limits_{i=0}^{k} (-1)^{i} \begin{bmatrix} m \\ k \end{bmatrix}_q \begin{bmatrix} k \\ i \end{bmatrix}_q q^{i(k+n-m)- \frac{i(i-1)}{2}}\, ch_{gr}W_{loc}(n-m+2(k-i))  \\
				= \sum\limits_{i=0}^{m} \sum\limits_{k=i}^{m} (-1)^{i} \begin{bmatrix} m \\ k \end{bmatrix}_q \begin{bmatrix} k \\ i \end{bmatrix}_q q^{i(k+n-m)- \frac{i(i-1)}{2}}\, ch_{gr} W_{loc}(n-m+2(k-i)) \\
				= \sum\limits_{i=0}^{m} \sum\limits_{k=0}^{i} (-1)^{k} \begin{bmatrix} m \\ m-i+k \end{bmatrix}_q \begin{bmatrix} m-i+k \\ k \end{bmatrix}_q q^{k(n-i+k)- \frac{k(k-1)}{2}}\, ch_{gr}W_{loc}(n+m-2i)\\
				= \sum\limits_{i=0}^{m} \begin{bmatrix} m \\ i \end{bmatrix}_q ( \sum\limits _{k=0}^{i} (-1)^{k}  \begin{bmatrix} i \\ k \end{bmatrix}_q q^{k(n-i+k)-\frac{k(k-1)}{2}}\,) ch_{gr}W_{loc}(n+m-2i)  	
			\end{array}$$
			
			To complete the proof of the lemma it suffices to show
			\begin{equation}\label{c identity} \sum_{k=0}^{i}  (-1)^{k}  \begin{bmatrix} i \\ k \end{bmatrix}_q q^{k(n-i+k)-\frac{k(k-1)}{2}} = \begin{bmatrix} n \\ i \end{bmatrix}_q (1-q)\dots (1-q^i ). \end{equation} 
			
			We prove \eqref{c identity} by using the following combinatorial identity:
			\begin{equation} \label{c identity 2}\sum\limits_{k=0}^{i} (-1)^{i-k} \begin{bmatrix} i \\ k \end{bmatrix}_q q^{{i-k\choose 2}} x^k= (x-1)(x-q) \cdots (x-q^{i-1})  
			\end{equation}
			Putting $x= q^n $ in \eqref{c identity 2}, we  get
			
			$$\begin{array}{rl}
				\sum\limits_{k=0}^{i} (-1)^{i-k} 
				\begin{bmatrix} i \\ k \end{bmatrix}_q q^{{i-k\choose 2}} q^{nk} 
				&= (q^n -1)(q^{n} -q)\cdots (q^n - q^{i-1}) \\ 
				\sum\limits_{k=0}^{i} (-1)^{i-k} \begin{bmatrix} i \\ k \end{bmatrix}_q 
				q^{(i-k)(i-k-1)/2} q^{nk} 
				&= (-1)^{i}q^{i(i-1)/2}(1-q^n )(1-q^{n-1}) \cdots (1-q^{n-i+1})  \\ 
				\sum\limits_{k=0}^{i}  (-1)^{k}  \begin{bmatrix} i \\ k \end{bmatrix}_q q^{k(n-i+k)- k(k-1)/2}  &=    (1-q^n ) \cdots (1-q^{n-i+1}) = \begin{bmatrix} n \\ i \end{bmatrix}_q (1-q)\dots (1-q^i ). \quad 
			\end{array}$$
			This establishes \eqref{c identity} and hence the lemma.
			\endproof
			
			\subsection{Ordered basis of the graded weight spaces of $W_{loc}(m)$} We have seen in \lemref{CV iso}(b), for $m\in \mathbb Z_+$, $W_{loc}(m)$ has a basis $\mathbb B(1^n)$ which is  indexed by the set $J(1^n)\cup \{\emptyset\}$. Define a function $|.|: J(1^n)\cup\{\emptyset\}\rightarrow \mathbb Z$  such that
			$$|\emptyset|=0, \qquad \quad |\boi| :=\sum\limits_{k=1}^m i_k, \quad \forall\, \boi =(i_1,\cdots,i_m)\in J(1^m).$$ Defining  in  $U(\mathfrak n^-[t])$ the elements  
			$$ y({{0}},\emptyset)=1, \qquad \quad 
			y(|\boi|,\boi):= (y_0)^{i_1}(y_1)^{i_2}\cdots(y_{m-1})^{i_m} ,\quad \forall\, \boi =(i_1,\cdots,i_m)\in J(1^m),$$ we see that, 
			$\mathbb B(1^n)=\{ y(|\boi|,\boi)w_m: \boi\in J(1^m)\cup\{\emptyset\}\}.$
			We define an ordering on $J(1^m)\cup \{\emptyset\}$ as follows. By fiat, $\emptyset >\boi$ for all $\boi\in J(1^m)$ and given $\boi,\boj\in J(1^m)$, we say, 
			$\boi>\boj$, if either $|\boi|< |\boj|$ or  $|\boi|=|\boj|$ 
			and there exists $1\leq k\leq m$ such that $i_k>j_k$ and   $i_s=j_s$ for $k+1\leq s\leq m$. This clearly induces an ordering in $\mathbb B(1^m)$, thereby making $\mathbb B(1^m)$ an ordered basis of $W_{loc}(m)$. \\
			\noindent On the other hand it was proved in \cite{MR2271991} that :
			\begin{proposition}\label{Wm.basis}  Given a positive integer $m$, if  
				$$F(m)=\{(k,\underline{s}): k\in \mathbb N, \underline{s}=(s_1,\cdots,s_k)\in {\mathbb Z}^k,0\leq s_i\leq m-k \text{ for } 1\leq i\leq k \}. $$ Then the set $\mathbb{B}(m)= \{y(k,\underline{s})w_{m}: (k,\underline{s}) \in F(m)\cup\{(0,\emptyset)\}\}$ is a basis of $W_{loc}(m)$, where $y(k,\underline{s}) = (y\otimes t^{s_1})\cdots(y\otimes t^{s_k})$ and $y(0,\emptyset)w_m=w_m$. \end{proposition}
			
			Thus using standard q-binomial theory, it follows from \propref{Wm.basis} 
			that the number of $l$-graded elements  of weight $m\omega-k\alpha$ in $\mathbb B(m)$ and hence in $\mathbb B(1^m)$ is equal to the coeffient of $q^l$ in the polynomial $\begin{bmatrix} m \\ k \end{bmatrix}_{q}$. 
			\label{ordered.basis.Wm}

			\subsection{Filtration of $W_{loc}(n)\otimes W_{loc}(m)$}
			\label{tensor char} 
			The following is the main result of this section.
			\begin{theorem}\label{tensor filtration}
				Let $n,m\in \mathbb N$ and $n \geq m$. The $\mathfrak{sl}_2[t]$-module $W_{loc}(n)\otimes W_{loc}(m)$ admits a filtration whose successive quotients are isomorphic to truncated local Weyl modules $$\tau_{k_r} W_{loc}([m+n-2r], n-r), \quad 0\leq r\leq m, \, 0\leq k_r\leq r(m-r).$$
			\end{theorem}
			\proof Given $n,m\in \mathbb N$ with $n\geq m$, let $w_n$ be the generator of $W_{loc}(n)$ and $\mathbb B(1^m)$ the ordered basis of $W_{loc}(m)$ as described in \secref{ordered.basis.Wm}. For each element $\boi\in J(1^m)\cup\{\emptyset\}$, 
			set 
			$$\begin{array}{c}v_{n,m}(k,\boi) = w_n\otimes y(k,\boi)w_m, \qquad \text{ where } k=|\boi|,\\ \\ \mathcal B(m,n)^{\boi\leq} = \{v_{n,m}(l,\boj): \boj\in J(1^m)\cup\{\emptyset\}, \boj\geq \boi\}
			\end{array}$$ and let $T(k,\boi)$ be the submodule of $W_{loc}(n)\otimes W_{loc}(m)$ generated by $\mathcal B(m,n)^{\boi\leq}$.\\

			\noindent Let $\boe_i[m]$ denote the $m$-tuple with 1 in the $i^{th}$-position and zero elsewhere. By \lemref{gen.W_loc.x.V(p)}, $$T(m,m\boe_1[m]) = W_{loc}(n)\otimes W_{loc}(m),$$ and from the definition of the submodules $T(|\boi|,\boi)$  it is clear that 
			$$T(m,m\boe_1[m])\supset T(m-1,(m-1)\boe_1[m])\supset \cdots \supset 
			T(1,\boe_{m-1}[m])\supset T(1,\boe_{m}[m])\supset T(0,\emptyset)\supset 0,$$ is a descending chain of submodules of $W_{loc}(n)\otimes W_{loc}(m)$. \\
			
			\noindent Additionally, since $x.y(k,\boi).w_m$ lies in the subspace of $W_{loc}(m)$ spanned by the vectors $\{y(k-1,\boj)w_m: \boj\in J(1^m)\}$ and 
			for $r>0$ , $\boi\in J(1^m)$ with $|\boi|=k$,
			$$(h\otimes t^r) y(k,\boi).w_m = \sum\limits_{\{1\leq j\leq m-r \,:\, 0<i_j<m-r\}} -2 y(k,\boi+(\boe_{i_j+r}[m]-\boe_{i_j}[m])).w_m,$$   using the fact that $x.w_n=0,$ and $h\otimes t^r.w_n=0$, for $r>0$, it follows that   
			\begin{equation}\begin{array}{c}
					x.v_{n,m}(1,\boi) =0 ,\qquad \quad 
					x.v_{n,m}(k,\boi) \in  T(k-1,(k-1)\boe_1[m])  \quad \text{ when } k>1,  \\ \\
					(h \otimes t^{r})v_{n,m}(k,\boi) \in \sum\limits_{\boj>\boi}\, T(k,\boi).\\
				\end{array}\label{xh.rel}\end{equation} 
			
			\noindent As $J(1^m)$, and hence $\mathbb B(1^m)$ is totally ordered, given 
			$\boi\in J(1^m)$ with $|\boi|=k$ there exists a unique pair $(\hat{k},\hat\boi)\in \mathbb Z\times J(1^m)$ such that $y(\hat{k},\hat{\boi})>y(k,\boi)$ and if 
			$y(l,\boj)\in \mathbb Z\times J(1^m)$ is such that $y(l,\boj)>y(k,\boi)$, then 
			$y(l,\boj)\geq y(\hat{k},\hat{\boi}).$ 
			In particular, as  $y(k-1,(k-1)\boe_1[m])>y(k,\boi)$, in view of \eqref{xh.rel}
			we see that,   
			$\dfrac{T(k,\boi)}{T(\hat{k},\hat{\boi})}$ is a highest weight module generated by the image of $v_{n,m}(k,\boi)$ under the natural surjection from $T(k,\boi)$ onto $\dfrac{T(k,\boi)}{T(\hat{k},\hat{\boi})}.$ \\
			\\ Since $n\geq m$, from the defining relations of a local Weyl module 
			we have $y_nv_{n,m}(k,\boi)=0$. Further,
			$$(h \otimes 1)v_{n,m}(k,\boi) = h.w_{n} \otimes y(k,\boi)w_{m} +  w_{n} \otimes h.y(k,\boi)w_{m} 
			= (n+m-2k) (w_{n} \otimes y(k,\boi)w_{m}).$$ 
			Hence, for all $\boi\in J(1^m)$,   the module  
			$\dfrac{T(k,\boi)}{T(\hat{k},\hat{\boi})}$ is a quotient of the truncated  
			Weyl module $\tau^\ast_{g_{\boi}}W_{loc}(m+n-2k,n)$, where 
			$g_{\boi}=\sum\limits_{k=1}^m (k-1)i_k$, $T(0,\emptyset)$ is a quotient of $W_{loc}(n+m,n)$ and 
			\begin{equation} \label{filter char} \, ch_{gr} W_{loc}(n)\otimes W_{loc}(m) = 
				ch_{gr} T(0,\emptyset) + \sum\limits_{\boi\in J(1^m)} q^{g_{\boi}}\, ch_{gr} \dfrac{T(|\boi|,\boi)}{T(|\hat{\boi}|, \hat{\boi})}.\end{equation}
			
			\noindent Now comparing with the character formulas \eqref{local tensor char} 
			and \eqref{filter char}, 
			we see that $T(0,\emptyset)$ is isomorphic $V(2^{n},1^{m-n})$ which by \secref{truncated} is isomorphic to $W_{loc}(m+n,n)$.\\
			\\
			On the other hand for $\boi\in J(1^m)$, with $|\boi|=k$ and $g_\boi =l$,  as
			$y(|\boi|,\boi)w_m$ is a $l$-graded vector of weight $m\omega-k\alpha$, the quotient module,
			$\dfrac{T(|\boi|,\boi)}{T(|\hat{\boi}|, \hat{\boi})}$ is a quotient of $\tau_l W(m+n-2k,n)$. By \secref{ordered.basis.Wm}, the number of l-graded vectors of weight $m\omega-k\alpha$ in $\mathbb B(1^m)$ is equal to the coefficient of $q^k$ in the polynomial $\begin{bmatrix} m \\ k \end{bmatrix}_{q}$. Hence comparing the character formulas \eqref{local tensor char} 
			and \eqref{filter char}  and the grades of the highest weight generating vectors the modules $\dfrac{T(|\boi|,\boi)}{T(|\hat{\boi}|, \hat{\boi})}$,  we see that  for any $0\leq k\leq m$,
			$\dfrac{T(k,\boi)}{T(\hat{k},\hat{\boi})}$ is isomorphic to the CV module $V(2^{m-k}, 1^{n-m})$ which,  by \secref{truncated}, is a quotient of $W_{loc}(n+m-2k,n)$and is isomorphic to the truncated Weyl module $W_{loc}(n+m-2k,n-k)$. This completes the proof of the theorem. \endproof

			\subsubsection{} For $m,n, s\in \mathbb N$, such that $m+n-2s\geq 0$, define a polynomial in indeterminate $q$ by: $$[W_{loc}(n)\otimes W_{loc}(m) : D(2,m+n-2s)]_q = \sum\limits_{p\geq 0} [W_{loc}(n)\otimes W_{loc}(m) : \tau_p D(2,m+n-2s)]q^p, 
			$$ where $[W_{loc}(n)\otimes W_{loc}(m) : \tau_p D(l,m+n-2s)]$ is the multiplicity of $\tau_p D(l,m+n-2s)$ in $W_{loc}(n)\otimes W_{loc}(m)$. The polynomial $[W_{loc}(n)\otimes W_{loc}(m) : D(2,m+n-2s)]_q$  is called the graded multiplicity of 
			$D(2, m+n-2s)$ in level $2$ Demazure flag of  $W_{loc}(n)\otimes W_{loc}(m).$
			\begin{corollary} Given two positive integers $m,n$, the $\mathfrak{sl}_2[t]$-module $W_{loc}(n)\otimes W_{loc}(m)$ has a level 2 Demazure flag and 
				$$\begin{array}{l}[W_{loc}(n)\otimes W_{loc}(m) : D(2,m+n-2s)]_q  \\
					=\left\{\begin{array}{ll} \sum\limits_{k=0}^{\min\{s,m\}} q^{(s-k)(m-k+\lceil\frac{n-m}{2}\rceil)}\begin{bmatrix} m\\ k\end{bmatrix} \begin{bmatrix} \lfloor\frac{n-m}{2}\rfloor \\ s-k \end{bmatrix}_{q}, & 0\leq s\leq \lfloor \frac{n-m}{2}\rfloor,\\
						\sum\limits_{k=0}^{\min\{m-j,\lfloor\frac{n-m}{2}\rfloor\}} q^{(\lfloor\frac{n-m}{2}\rfloor-k)(m-k-j+\lceil\frac{n-m}{2}\rceil)}
						\begin{bmatrix} m\\ k+j\end{bmatrix} 
						\begin{bmatrix} \lfloor\frac{n-m}{2}\rfloor \\ \lfloor\frac{n-m}{2}\rfloor-k \end{bmatrix}_{q}, & s=j+\lfloor \frac{n-m}{2}\rfloor, \,  1\leq j\leq m
					\end{array}\right.\end{array} $$ for $m+n-2s\geq 0$.
			\end{corollary}
			
			\proof By \thmref{tensor filtration}, the tensor product product module $W_{loc}(n)\otimes W_{loc}(m)$ has a filtration by submodules whose 
			successive quotients are isomorphic to translates of modules of the form $V(2^a,1^b)$. Hence by \thmref{Theorem 3.3}, $W_{loc}(n\omega)\otimes W_{loc}(m\omega)$ has a level 2 Demazure flag. Now using \lemref{Blanton.result} and \lemref{char V(a,b)}, we see that
			$$\begin{array}{ll}
				ch_{gr}W_{loc}(n){\otimes}W_{loc}(m)&= \sum\limits_{k=0}^{m} 
				\begin{bmatrix} m \\ k \end{bmatrix}_{q} ch_{gr}\, V(2^{m-k},1^{n-m})\\
				&=\sum\limits_{k=0}^m \begin{bmatrix} m \\ k \end{bmatrix}_{q} 
				\left(\sum\limits_{r=0}^{\lfloor\frac{n-m}{2}\rfloor} q^{r(m-k+\lceil\frac{n-m}{2}\rceil)} 
				\begin{bmatrix} \lfloor\frac{n-m}{2}\rfloor \\ r \end{bmatrix}_{q} ch_{gr} D(2,n+m-2(k+r)) \right)\\
			\end{array}$$
			Hence, we have
			$$\begin{array}{ll}[W_{loc}(n\omega)\otimes W_{loc}(m\omega): D(2,m+n-2s)]_q  
				&= \sum\limits_{k=0}^{\min\{s,m\}} q^{(s-k)(m-k+\lceil\frac{n-m}{2}\rceil)}
				\begin{bmatrix} m\\ k\end{bmatrix} 
				\begin{bmatrix} \lfloor\frac{n-m}{2}\rfloor \\ s-k \end{bmatrix}_{q}, 
			\end{array} 
			$$  if $s\in \mathbb N$ is such that $0\leq s\leq \lfloor \frac{n-m}{2}\rfloor$  and 
			$$\begin{array}{l}[W_{loc}(n\omega) \otimes W_{loc}(m\omega): D(2,m+n-2s)]_q  \\
				=\sum\limits_{k=0}^{\min\{m-j,\lfloor\frac{n-m}{2}\rfloor\}} q^{(\lfloor\frac{n-m}{2}\rfloor-k)(m-k-j+\lceil\frac{n-m}{2}\rceil)}
				\begin{bmatrix} m\\ k+j\end{bmatrix} 
				\begin{bmatrix} \lfloor\frac{n-m}{2}\rfloor \\ \lfloor\frac{n-m}{2}\rfloor-k \end{bmatrix}_{q},
			\end{array} $$ 
			when $s=j+\lfloor \frac{n-m}{2}\rfloor$ for some $1\leq j\leq m$. \endproof
			
			In \cite{9} the graded character formula for $D(2, 2n+r) \otimes W_{loc}(m)$ has been obtained when $n\geq m$. We can extend the methods employed in proving Theorem \ref{tensor filtration} to show that the tensor product of a level 2 Demazure module with a local Weyl module has a decreasing chain of submodules, with each successive quotient isomorphic to a CV-module of the form $V(3^a, 2^b, 1^c)$. Consequently, it follows from \propref{Theorem 3.3} that the module $D(2, 2n + r) \otimes W_{\text{loc}}(m)$ admits a level 3-Demazure flag.

			\section{Appendix}
			In the appendix we give a proof of Lemma \ref{tensor of two local} using Pieri formulas (\cite[Chapter VI]{MR3443860}). \\
			
			\noindent Let $t, q$ be independent indeterminates and  $\mathbb C(q, t)$ be the field of rational functions in $q$ and $t$. Let $\{x_i: 1\leq i\leq n+1\}$ be a set of indeterminates and $\mathbb C(q,t)[x_1,\cdots,x_{n+1}]$ be the ring of polynomials in   $\{x_j: 1\leq j\leq n+1\}$ with coefficients in $\mathbb C(q,t)$. Let $\boldsymbol{x} = (x_1 , x_2 , \dots , x_{n+1})$ be a basis for the ring of symmetric polynomials in $\mathbb{C}(q,t)[x_1 , x_2 , \dots , x_{n+1}]$.\\
			\\
			For $\lambda\in P^+$, a family of orthogonal symmetric polynomials $P_{\lambda}(\boldsymbol{x}; q, t)$ was introduced  in \cite{MR3443860}, and it was shown in \cite{MR1771615} that the character of a level one Demazure module of highest weight $\lambda$ is given by the specialized Macdonald polynomial $P_{\lambda}(\boldsymbol{x}; q, 0)$. On the other hand it was proved in  \cite{MR2323538} that
			a local Weyl module for a simply-laced current Lie algebra is isomorphic to a Demazure module of level one. Thus in case of $\mathfrak{sl}_{2}[t]$, we see that for $m\in \mathbb Z_+$,  $ch_{gr} W_{loc} (m) = P_{(m,0)}(\boldsymbol{x};q,0)$. We use this fact to give an alternate proof of \lemref{tensor of two local}.\\
			
			Set  \begin{align*}
				(q;q)_{n} &= (1-q)(1-q^2 ) \dots (1-q^n ) \quad  \\
				(t;q)_n &= (1-t)(1-tq) \dots (1-tq^{n-1}) \quad \text{for } |q|<1 \\
				(t;q)_{\infty} &= \prod_{i=0}^{\infty}(1-tq^{i})
			\end{align*}
			By \cite[Chapter VI, Equation 4.9]{MR3443860}, 
			\begin{equation}\label{gm equation}
				P_{(m)}(\boldsymbol{x};q,t)= \frac{(q;q)_{m}}{(t;q)_{m}}g_{m}(\boldsymbol{x},q,t) 
			\end{equation}
			where $g_{m}(\boldsymbol{x};q,t)$ denotes the coefficient of $y^m $ in the power-series expansion of the infinite product 
			$$\prod_{i\geq 1} \frac{(tx_{i} y;q)_{\infty}}{(x_{i}y;q)_{\infty}} = \sum_{m \geq 0}g_{m}(\boldsymbol{x};q,t)y^{m}$$
			
			\noindent Putting $t=0$ in \eqref{gm equation}, we  get
			\begin{equation} \label{P in terms of g}
				P_{(m)}(\boldsymbol{x};q,0)= (q;q)_{m}g_{m}(\boldsymbol{x},q,0)
			\end{equation}
			
			\vspace{.15cm}
			
			\noindent We need the following result on Pieri formula for our proof of \lemref{tensor of two local}.
			\begin{lemma}  \cite[Theorem 6.24]{MR3443860} 
				Let $m,n \in \mathbb Z_{+}$, with $n \geq m$. Given a partition $\lambda$ of $m+n$ containing the partition $(n)$ and having atmost two parts, we have,
				$$P_{n}(\boldsymbol{x};q,0) g_{m}(\boldsymbol{x};q,0) = \sum_{\lambda - (n) = (m)} \phi_{\lambda/ (n)} P_{\lambda}(\boldsymbol{x};q,0)$$ 
				and the coefficients are given by 
				$$ \phi_{\lambda / (n)} = \prod_{s \in C_{\lambda / (n)}} \frac{b_{\lambda}(s)}{b_{(n)}(s)}$$ 
				where  $C_{\lambda /(n)}$ denote the union of columns that intersect $\lambda-(n)$, and 
				$$ b_{\lambda}(s) = b_{\lambda}(s;q,t) = \begin{cases} 
					\frac{1-q^{a_{\lambda}(s)}t^{l_{\lambda}(s)+1}}{1- q^{a_{\lambda}(s)+1}t^{l_{\lambda}(s)}} & \quad \text{if } s \in \lambda, \\
					1 & \quad \text {otherwise} 
				\end{cases} $$
			\end{lemma}
			\subsection*{Proof of \lemref{tensor of two local}.}
			Since for $m\in \mathbb Z_+$, $ch_{gr} W_{loc} (m) = P_{(m,0)}(\boldsymbol{x};q,0)$, we have
			\begin{align*}
				ch_{gr}(W_{loc}(n) \otimes W_{loc}(m)) &= ch_{gr}W_{loc}(n)ch_{gr}W_{loc}(m)\\
				& = P_{(n,0)}(\boldsymbol{x};q,0)P_{(m,0)}(\boldsymbol{x};q,0)\\
				& = (q;q)_{m}P_{(n,0)}(\boldsymbol{x};q,0)g_{m}(\boldsymbol{x};q,0) \quad(\text{Using \ref{P in terms of g}})\\
				& = (q;q)_{m} \sum_{\lambda - (n) = (m)} \phi_{\lambda/ (n)}(q,0) P_{\lambda}(\boldsymbol{x};q,0)\\
				& = (q;q)_{m} \sum_{\lambda - (n) = (m)} \prod_{s \in C_{\lambda / (n)}} \frac{b_{\lambda}(s;q,0)}{b_{(n)}(s;q,0)} P_{\lambda}(\boldsymbol{x};q,0)
					\end{align*}
				\begin{align*}
				& = (q;q)_{m} (\frac{1}{(1-q)(1-q^2 ) \dots (1-q^m )} P_{(n+m,0)}(\boldsymbol{x};q,0) \\
				& + \frac{1-q^n }{(1-q)(1-q^2 ) \dots (1-q^{m-1} )(1-q)} P_{(n+m-1,1)}(\boldsymbol{x};q,0) \\
				&+ \frac{(1-q^n)(1-q^{n-1}) }{(1-q)(1-q^2 ) \dots (1-q^{m-2} )(1-q)(1-q^2 )} P_{(n+m-2,2)}(\boldsymbol{x};q,0)\\
				&+ \dots + \frac{(1-q^n )(1-q^{n-1}) \dots (1-q^{n-m+1})}{(1-q)(1-q^2 ) \dots (1-q^m )}
				P_{(n-m,m)}(\boldsymbol{x};q,0))\\
				& = P_{(n+m,0)}(\boldsymbol{x};q,0) + \frac{(1-q^n)(1-q^m )}{(1-q)} P_{(n+m-1,1)}(\boldsymbol{x};q,0) \\
				& + \frac{(1-q^n)(1-q^{n-1})(1-q^m )(1-q^{m-1}) }{(1-q)(1-q^2 )} P_{(n+m-2,2)}(\boldsymbol{x};q,0)\\
				& + \dots + \frac{(1-q^n )(1-q^{n-1}) \dots (1-q^{n-m+1})(1-q)(1-q^2 ) \dots (1-q^m )}{(1-q)(1-q^2 ) \dots (1-q^m )}\\
				& P_{(n-m,m)}(\boldsymbol{x};q,0)\\
			\end{align*} Hence,
			\begin{align*}
				ch_{gr}(W_{loc}(n) \otimes W_{loc}(m))  & = \sum_{k=0}^{m} \begin{bmatrix} n\\ k \end{bmatrix}_{q} \begin{bmatrix} m\\ k \end{bmatrix}_{q} 
				(1-q)(1-q^2 ) \dots (1-q^k ) P_{(n+m-k,k)}(\boldsymbol{x};q,0)\\
				& = \sum_{k=0}^{m} \begin{bmatrix} n \\ k \end{bmatrix}_q \begin{bmatrix} m \\ k \end{bmatrix}_q (1-q)\dots (1-q^k ) ch_{gr} W_{loc}(n+m-2k). \qquad \qquad  \qed
			\end{align*}
			
			\section*{\textbf{Acknowledgements}}
			The first author would like to thank Professor S.Vishwanath and Professor Vyjayanthi Chari for useful discussions. The first author is grateful to IISER Mohali for PhD fellowship and support.

			
		\end{document}